\documentstyle[12pt, righttag]
{amsart}
\input amssym.def
\input amssym

\newcommand{\wt}{\mbox{\normalshape wt}}
\newcommand{\spa}{\mbox{\normalshape span}}
\newcommand{\Res}{\mbox{\normalshape Res}}
\newcommand{\End}{\mbox{\normalshape End}}
\newcommand{\Ind}{\mbox{\normalshape Ind}}
\newcommand{\Hom}{\mbox{\normalshape Hom}}
\newcommand{\Mod}{\mbox{\normalshape Mod}}
\newcommand{\m}{\mbox{\normalshape mod}\ }

\def \Z{\Bbb Z}
\def \M{\Bbb M}
\def \C{\Bbb C}

\def \Q{\Bbb Q}

\def \<{\langle}
\def \o{\omega}
\def \O{\Omega}
\def \M{{\cal M}}
\def \1t{\frac{1}{T}}
\def \>{\rangle}

\def \l{\lambda }
\def \L{\Lambda }

\def \o{\omega }

\def \v{vertex operator algebra\ }

\def \be{\begin{equation}\label}
\def \ee{\end{equation}}
\def \qed{\mbox{ $\square$}}
\def \pf {\noindent {\bf Proof:} \,}
\def \bl{\begin{lem}\label}
\def \el{\end{lem}}
\def \ba{\begin{array}}
\def \ea{\end{array}}
\def \bt{\begin{thm}\label}
\def \et{\end{thm}}

\def \br{\begin{rem}\label}
\def \er{\end{rem}}
\def \ed{\end{de}}
\def \bp{\begin{prop}\label}
\def \ep{\end{prop}}

\newtheorem{thm}{Theorem}[section]
\newtheorem{prop}[thm]{Proposition}

\newtheorem{lem}[thm]{Lemma}
\newtheorem{rem}[thm]{Remark}
\newtheorem{de}[thm]{Definition}

\numberwithin{equation}{section}

\begin{document}
\title[Twisted representations of vertex operator algebras
]{Twisted representations of vertex operator algebras}
\author{Chongying Dong, Haisheng Li and Geoffrey Mason}
\address{Department of Mathematics,\,University of California,\,Santa Cruz}
    \email{dong@@cats.ucsc.edu\\
    hli@@cats.ucsc.edu\\
    gem@@cats.ucsc.edu}
\thanks{C.D is partially supported by NSF grant DMS-9303374 and a
research grant from the Committee on Research, UC Santa Cruz.}
\thanks{G.M. is partially supported by NSF grant
DMS-9401272 and a research grant from the Committee on Research, UC
Santa Cruz.}
\bibliographystyle{alpha}
	\maketitle

\begin{abstract} Let $V$ be a vertex operator algebra and $g$
an automorphism of finite order. We construct an associative
algebra $A_g(V)$ and a pair of functors between
the category of $A_g(V)$-modules and a certain category
of admissible $g$-twisted $V$-modules. In particular, these
functors exhibit a bijection between the simple modules in each
category. We give various applications, including the fact that
the complete reducibility of admissible $g$-twisted modules implies
both the finite-dimensionality of homogeneous spaces and the finiteness
of the number of simple $g$-twisted modules.
\end{abstract}

\section{Introduction}

The study of vertex operator algebras has come to play a significant
role in disparate areas such as conformal field theory [MS], Moonshine
and the Monster [FLM], [B] and elliptic cohomology [T]. Although fairly
recent in origin, some of the main problems in the theory of vertex
operator algebras are quite classical in nature and concern representation
theory.

Let $V$ be a \v and $G$ be a finite automorphism group of $V.$ Then the
space of
$G$-invariants $V^G$ is itself
a vertex operator algebra.  It is natural
to try to understand
various module categories for $V^G.$  This is so-called
{\em orbifold theory} in the physical literature [DHVW], [DVVV].
One of the main new features of
orbifold theory is the introduction of {\em twisted modules} or
{\em twisted sectors.} Essentially, these are spaces which admit
vertex operators indexed by elements of $V$ and satisfying analogues of
the Jacobi identity which are ``twisted'' by elements $g\in G.$ Moreover
they restrict to ``ordinary'' modules for $ V^G.$

Although, for certain orbifolds, some success
has been achieved in the study of these objects
in the physics literature [DVVV], [DGM], the mathematical
investigation of abstract orbifold models has been hampered
by a lack of understanding of
the theory of twisted representations of vertex operator algebras.
The goals of the present paper are to alleviate this situation.

Given a \v $V$ and  automorphism $g$ of finite order
$T,$ we will construct an associative
algebra $A_g(V)$ with the property that there is a bijective correspondence
between simple $A_g(V)$-modules and simple {\em admissible} $g$-twisted
$V$-modules. These latter objects are twisted analogues of Zhu's definition
of $V$-modules [Z]. They carry a grading by $\1t\Z_+,$ but
the homogeneous spaces are neither assumed to be of finite dimension, nor
induced from the eigenvalues of the $L(0)$ operator. If these latter
conditions hold, then we have an ({\em ordinary}) $g$-twisted module
as defined in [FFR] and [D1].
In fact, the main concern of the paper is the
construction of a pair
of functors $L,$ $\O$ which we display follows:
$$A_g(V)-\Mod
\begin{array}{l}
\stackrel{L}{\longrightarrow} \\
\stackrel{\Omega}{\longleftarrow}
\end{array}
\mbox{Adm}-g-V-\Mod$$
Thus the functor $L$ constructs a certain
admissible $g$-twisted $V$-module $L(U)$ from a given $A_g(V)$-module
$U,$ whilst the functor $\Omega$ does the opposite.
Moreover we have $\O\circ L\cong id.$ Furthermore
$L$ and $\O$ induce bijections on the simple objects of each category.
Because of the failure of complete reducibility of appropriate modules
one cannot expect $\O$ and $L$ to be mutually inverse
categorical equivalences in general, though we are able
to prove this (Theorem \ref{t7.2}) for the full subcategories of completely
reducible objects.

There is an important application of our theory to $g$-{\em rational}
vertex operator algebras; these are vertex operator algebras
such that every  admissible $g$-twisted module is completely
reducible. We show (Theorem \ref{t8.1}) that such vertex operator algebras
necessarily
have only finitely many inequivalent
simple admissible $g$-twisted modules, and that
every such module is an ordinary $g$-twisted module. So for $g$-rational
vertex operator algebras, $L$ and $\O$ induce mutually inverse
categorical equivalences between the categories of finitely
generated $A_g(V)$-modules and ordinary $g$-twisted $V$-modules.

We have already alluded to Zhu's work [Z]. Our theory includes that of
Zhu if we take $g=1,$ but
 even in this case our work leads to a strengthening of some of his results
as well as a simplification in the proofs. One of our main ideas, which
goes back to [L2] if $g=1,$ is the introduction
of a certain Lie algebra $V[g]$ into the proceedings. This allows
us to replace Zhu's use of correlation functions, which
is quite difficult,  with more familiar methods of Lie theory
(induced modules, PBW theorem).

We have already made use
of our results in several papers [DLM1]-[DLM2], [DM1]-[DM2],
and expect that the study of $g$-twisted modules will
lead to a proof of the generalized Moonshine conjectures when applied to the
action of the Monster on the Moonshine Module.

The paper is organized as follows: in Section 2 we introduce the
algebra $A_g(V).$ In Section 3 we discuss the various kinds of
$g$-twisted $V$-modules that we need to deal with. In Section 4 we
construct the Lie algebra $V[g]$ and show that a {\em weak}
$g$-twisted $V$-module is a $V[g]$-module. Then in Section 5 we
construct the functor $\O;$ it is obtained essentially as the space of
lowest weight vectors for $V[g].$ Section 6 is technically the most
difficult. We construct the functor $L,$ which entails the
construction of a certain graded $V[g]$-module $L(U)$ from a given
$A_g(V)$-module $U$ and then verifying the twisted Jacobi identity.
This is never easy! We also construct (Theorems \ref{t6.1} and
\ref{t6.3}) a certain universal object $\bar M(U)$ in the category of
admissible $g$-twisted $V$-modules, and which has
$L(U)$ as a quotient. Thus $\bar M(U)$ is a sort of ``generalized''
Verma module. In Section 7 we prove that $L$ and $\O$ are equivalences
when restricted to the subcategory of completely reducible
objects. Section 8 is concerned with $g$-rational vertex operator
algebras and includes the results already mentioned.
Section 9 contains some useful
applications. It should be emphasized that it remains a conjecture
that non-zero $g$-twisted $V$-modules always exists; we prove that this is
so if $A_g(V)$ is of finite dimension (Theorem \ref{t9.1}). We also
give some sufficient conditions for the complete reducibility of (admissible
and ordinary) $V$-modules.

We expect the reader to be familiar with the elementary theory of vertex
operator algebras as found, for example, in [FLM] and [FHL].

\section{The associative algebra $A_{g}(V)$}

We fix some notation which will be in force throughout the paper.
$(V,Y,{\bold 1},\omega)$ denotes, as usual, a
vertex operator algebra (cf. [B],
[FHL] and [FLM]) and $g$ is an automorphism of $V$ of finite order
$T.$  Denote the decomposition of $V$ into eigenspaces with respect to
the action of $g$ as
\begin{equation}\label{g2.1}
V=\oplus_{r\in \Z/T\Z}V^r
\end{equation}
where $V^r=\{v\in V|gv=e^{2\pi ir/T}v\}$.
(We habitually use $ r$ to denote both
an integer between $0$ and $T-1$ and its residue class \m $T$ in this
situation.)

We are going to construct an associative algebra $A_g(V)$ along the line of
Zhu's construction of his algebra $A(V)$ [Z]. Indeed if $g=1$ our algebra
$A_1(V)$ is {\em precisely} $A(V).$ In general one may consider Zhu's
algebra $A(V^0)$ associated with the vertex operator subalgebra $V^0$ of
$g$-invariants; our algebra $A_g(V)$ will be a certain quotient
of $A(V^0).$

For homogeneous $u\in V^r$ and $v\in V$ we  define
\begin{eqnarray}\label{g2.2}
u\circ_g v=\Res_{z}\frac{(1+z)^{{\wt u}-1+\delta_{r}+{r\over
T}}}{z^{1+\delta_{r}}}Y(u,z)v
\end{eqnarray}
where $\delta_{r}=1$ if $r=0$ and
$\delta_{r}=0$ if $r\ne 0$.

Let $O_g(V)$ be the linear span of all $u\circ_g v$ and define the linear space
$A_g(V)$ to be the quotient $V/O_{g}(V).$ We will usually write $A(V), O(V),
u\circ v,$ etc. when $g=1.$

\begin{lem}\label{l2.1} If $r\ne 0$ then $V^r\subseteq O_{g}(V).$
\end{lem}

\pf It suffices to show that $u\in O_g(V)$ whenever
$u\in V^r$ is homogeneous. In this case, take $v={\bold 1}$ in (\ref{g2.2})
to see that
$$u\circ_g {\bold 1}=\Res_{z}\frac{(1+z)^{{\wt u}-1+{r\over T}}}{z}Y(u,z){\bold
1}=u\in O_g(V).$$
The lemma follows. \qed

If we set $I=O_g(V)\cap V^0,$ it follows from Lemma \ref{l2.1} that $A_g(V)
\simeq V^0/I$ (as linear spaces). Notice that $O(V^0)\subset I,$ so that
$A_g(V)$ is a quotient of $A(V^0).$

Define a second product $*_g$ on $V$ as follows: with $r,u$ and $v$ as
above, set
\begin{equation}\label{a5.1}
u*_gv=\left\{
\begin{array}{ll}
\Res_z(Y(u,z)\frac{(z+1)^{{\wt}\,u}}{z}v)
 & {\rm if}\ r=0\\
0  & {\rm if}\ r>0.
\end{array}\right.
\end{equation}
Extend linearly to obtain a bilinear product  on $V$ which coincides with that
of Zhu (loc.cit.) on $V^0.$ We denote the product (\ref{a5.1})
by $u*v$ in this case. In this way $V$ becomes a (non-associative) algebra
with respect to $*_g.$ Note that if $u\in V^0$ then (\ref{a5.1}) may be
written in the form
\begin{equation}\label{g2.4}
u*_gv=\sum_{i=0}^{\infty}{\wt u\choose i}u_{i-1}v.
\end{equation}

Following Lemmas 2.1.2 and 2.1.3 of [Z], we get the following.
\begin{lem}\label{l2.2} (i) Assume that $u\in V^{r}$ homogeneous,
$v\in V$ and $m\ge n\ge 0.$ Then
$$\Res_{z}\frac{(1+z)^{{\wt}u-1+\delta_{r}+{r\over
T}+n}}{z^{m+\delta_{r}+1}}Y(u,z)v\in O_{g}(V).$$

(ii) Assume that $u,v\in V^0$ are homogeneous. Then
$$u*v-\Res_z\frac{(1+z)^{\wt v-1}}{z}Y(v,z)u\in O(V^0)$$
and

(iii) $u*v-v*u-\Res_z(1+z)^{\wt u-1}Y(u,z)v\in O(V^0).$
 \end{lem}

\begin{prop}\label{p2.7} (i) $O_{g}(V)$ is an two-sided
ideal of $V$ with respect to the product $*_g.$

(ii) If $I=O_g(V)\cap V^0$ then $I/O(V^0)$ is  a
two-sided ideal of $A(V^{0})$.
\end{prop}

\pf Notice that if $r\ne 0$ then $V^r*_gV=0$ by definition (\ref{a5.1}).
Similarly (\ref{a5.1}) shows that $V^0*_gV^r\subset V^r.$ Since
$O_g(V)=I\oplus(\oplus_{r=1}^{T-1}V^r)$ by Lemma \ref{l2.1}, we see that
parts (i) and (ii) are equivalent to each other and to the assertion that
$I$ is a 2-sided ideal of $V^0$ with respect to $*.$ We prove this latter
assertion.

Choose $c\in V^0$ homogeneous and $u\in I.$ We must show that $I$
contains both
\begin{eqnarray}
c*u=\Res_{z}\left(\frac{(1+z)^{{\wt}c}}{z}Y(c,z)u\right)
\end{eqnarray}
and
\begin{eqnarray}
u*c\equiv\Res_{z}\frac{(1+z)^{{\wt}c-1}}{z}Y(c,z)u\ (\m I).
\end{eqnarray}
(For the latter congruence use Lemma \ref{l2.2} (ii).) From (\ref{g2.2})
it suffices to take $u=a\circ_g b$ where $a\in V^r$ and
$b\in V^{T-r}$ are both homogeneous.
Set $x_{0}=c*u,$ $x_1=u*c$ and recall the Jacobi identity on $V:$
\begin{equation}\label{g2.7}
\begin{array}{c}
\displaystyle{z^{-1}_0\delta\left(\frac{z_1-z_2}{z_0}\right)
Y(c,z_1)Y(a,z_2)b-z^{-1}_0\delta\left(\frac{z_2-z_1}{-z_0}\right)
Y(a,z_2)Y(c,z_1)}b\\
\displaystyle{=z_2^{-1}\delta\left(\frac{z_1-z_0}{z_2}\right)
Y(Y(c,z_0)v,z_2)b}.
\end{array}
\end{equation}
Using (\ref{g2.7}) we have for $\varepsilon=0$ or $1:$
\begin{eqnarray*}
& &x_{\varepsilon}=\Res_{z_{1}}\frac{(1+z_{1})^{{\wt}c-\varepsilon}}{z_{1}}
Y(c,z_{1})\Res_{z_2}\frac{(1+z_2)^{{\wt}a-1+\delta_{r}+{r\over
T}}}{z_2^{1+\delta_{r}}}Y(a,z_2)b
\nonumber\\
& &=\Res_{z_{1}}\Res_{z_2}\frac{(1+z_{1})^{{\wt}c-\varepsilon}}{z_{1}}
Y(c,z_{1})\frac{(1+z_2)^{{\wt}a-1+\delta_{r}+{r\over
T}}}{z_2^{1+\delta_{r}}}Y(a,z_2)b\\
& &=\Res_{z_{1}}\Res_{z_2}\frac{(1+z_{1})^{{\wt}c-\varepsilon}}
{z_{1}}\frac{(1+z_2)^{{\wt}a-1+\delta_{r}+{r\over
T}}}{z_2^{1+\delta_{r}}}Y(a,z_2)Y(c,z_{1})b\nonumber\\
& & \ \ \ +\Res_{z_{1}}\Res_{z_2}\Res_{z_{0}}
\frac{(1+z_{1})^{{\wt}c-\varepsilon}}{z_{1}}\frac{(1+z_2)^{{\wt}a-1+\delta_{r}+{r\over T}}}{z_2^{1+\delta_{r}}}
z_2^{-1}\delta\left(\frac{z_{1}-z_{0}}{z_2}\right)Y(Y(c,z_{0})a,z_2)b\nonumber\\
& &=\Res_{z_{1}}\Res_{z_2}\frac{(1+z_{1})^{{\wt}c-\varepsilon}}
{z_{1}}\frac{(1+z_2)^{{\wt}a-1+\delta_{r}+{r\over
T}}}{z_2^{1+\delta_{r}}}Y(a,z_2)Y(c,z_{1})b\nonumber\\
& &\ \ \ +\Res_{z_2}\Res_{z_{0}}
\frac{(1+z_2+z_{0})^{{\wt}c-\varepsilon}}{z_2+z_{0}}\frac{(1+z_2)^{{\wt}a-1+\delta_{r}+{r\over T}}}{z_2^{1+\delta_{r}}}Y(Y(c,z_{0})a,z_2)b\nonumber\\
& &=\Res_{z_2}\frac{(1+z_2)^{{\wt}a-1+\delta_{r}+{r\over
T}}}{z_2^{1+\delta_{r}}}Y(a,z_2)
\Res_{z_{1}}\frac{(1+z_{1})^{{\wt}c-\varepsilon}}
{z_{1}}Y(c,z_{1})b\nonumber\\
& &\ \ \
+\sum_{i,j=0}^{\infty}(-1)^{j}\left(\begin{array}{c}{\wt}c-\varepsilon\\i\end{array}\right)\Res_{z_2}
\frac{(1+z_2)^{{\wt}a-1+\delta_{r}+{r\over T}+{\wt}c-\varepsilon -i}}
{z_2^{j+2+\delta_{r}}}Y(c_{i+j}a,z_2)b\nonumber\\
& &=\Res_{z_2}\frac{(1+z_2)^{{\wt}a-1+\delta_{r}+{r\over
T}}}{z_2^{1+\delta_{r}}}Y(a,z_2)
\Res_{z_{1}}\frac{(1+z_{1})^{{\wt}c-\varepsilon}}
{z_{1}}Y(c,z_{1})b\nonumber\\
& &\ \ \
+\sum_{i,j=0}^{\infty}(-1)^{j}\left(\begin{array}{c}{\wt}c-\varepsilon\\i\end{array}\right)\Res_{z_2}
\frac{(1+z_2)^{{\wt}(c_{i+j}a)-1+\delta_{r}+{r\over T}+j+1-\varepsilon}}
{z_2^{j+2+\delta_{r}}}Y(c_{i+j}a,z_2)b.\nonumber
\end{eqnarray*}
The resulting element is in $I$ by the definition of $O_g(V)$ and
Lemma \ref{l2.2} (i). The proof is complete. \qed

Our first main result is the following.
\begin{thm}\label{t2.4}  (i) The product $*_g$ induces the structure of an
associative algebra  on $A_{g}(V).$

(ii) The linear map
$$\phi:  a\mapsto e^{L(1)}(-1)^{L(0)}a$$
induces an anti-isomorphism $A_{g}(V)\to A_{g^{-1}}(V)$.

(iii) There are identities
$${\bold 1}*_g x\equiv x*_g{\bold 1}\equiv x\ (\m O_g(V))$$
$$\omega *_g x\equiv x*_g\omega\ (\m O_g(V))$$
for $x\in V.$
\end{thm}

\begin{rem} We have not ruled out the possibility that $A_g(V)$ is equal to
0, indeed this is a subtle
point related to the existence/non-existence of $g$-twisted sectors, as we
shall see below. Part (iii) says that the vacuum ${\bold 1}$ maps onto the
identity of $A_g(V)$ as long as $A_g(V)\ne 0.$ Similarly the image of the
Virasoro vector $\omega $ lies in the center of $A_g(V).$
\er

{\bf Proof of Theorem \ref{t2.4}:} Part (i) follows immediately from
Proposition \ref{p2.7} and Zhu's results [Z] that $A(V^0)$  is an associative
algebra with respect to $*.$ Similarly part (iii) follows from
Theorem 2.1.1 (2), (3) of [Z].

Since part (ii) follows for $g=1$ from [Z] (cf. equation (2.1.9) of that
paper, though no proof is given) the main point is to show that
$\phi$  maps
$O_{g}(V)\cap V^0$ into $O_{g^{-1}}(V)\cap V^0.$  First
recall the following conjugation formulas from [FHL]:
$$z^{L(0)}Y(a,z_{0})z^{-L(0)}=Y(z^{L(0)}a,zz_{0}),$$
$$e^{zL(1)}Y(a,z_{0})e^{-zL(1)}=Y\left(e^{z(1-zz_{0})L(1)}(1-zz_{0})^{-2L(0)}a,{z\over 1-zz_{0}}\right).$$
Then for homogeneous $a\in V^r$ and $b\in V^{T-r}$
we have:
\begin{eqnarray*}
& &\ \ \ \phi(a\circ_g b)\\
& &=\phi\left(\Res_{z}\frac{(1+z)^{{\wt}a-1+\delta_{r}+{r\over
T}}}{z^{1+\delta_{r}}}Y(a,z)b\right)
\nonumber\\
& &=e^{L(1)}(-1)^{L(0)}\Res_{z}\frac{(1+z)^{{\wt}a-1+\delta_{r}+{r\over
T}}}{z^{1+\delta_{r}}}Y(a,z)b
\nonumber\\
& &=\Res_{z}\frac{(1+z)^{{\wt}a-1+\delta_{r}+{r\over
T}}}{z^{1+\delta_{r}}}e^{L(1)}
Y((-1)^{L(0)}a,-z)(-1)^{L(0)}b\nonumber\\
& &=\Res_{z}\frac{(1+z)^{{\wt}a-1+\delta_{r}+{r\over T}}}{z^{1+\delta_{r}}}
Y\left(e^{(1+z)L(1)}(1+z)^{-2L(0)}(-1)^{L(0)}a,{-z\over
1+z}\right)e^{L(1)}(-1)^{L(0)}b.\nonumber\\
& &\mbox{}
\end{eqnarray*}
Replacing  $z$ with $\displaystyle{-{z_{0}\over 1+z_{0}}}$ and
using the residue formula for the change of variable (see [Z])
$$\Res_zg(z)=\Res_{z_0}(g(f(z_0)){d\over{dz_0}}f(z_0))$$
gives
\begin{eqnarray*}
& &\ \ \  \phi(a\circ_g b)\nonumber\\
& &=-\Res_{z_{0}}(1+z_{0})^{-{\wt}a+1-\delta_{r}-{r\over T}}
\left(\frac{-z_{0}}{1+z_{0}}\right)^{-1-\delta_{r}}
\frac{1}{(1+z_{0})^{2}}\cdot\nonumber\\
& &\ \ \ \cdot
Y\left(e^{(1+z_{0})^{-1}L(1)}(1+z_{0})^{2L(0)}(-1)^{L(0)}a,z_{0}\right)
 e^{L(1)}(-1)^{L(0)}b\nonumber\\
& &=(-1)^{{\wt}a+\delta_{r}}\Res_{z_{0}}\frac{(1+z_{0})^{{\wt}a-{r\over T}}}
{z_{0}^{1+\delta_{r}}}Y(e^{(1+z_{0})^{-1}L(1)}a,z_{0})e^{L(1)}(-1)^{L(0)}b\nonumber\\
& &=\sum_{j=0}^{\infty}{1\over j!}(-1)^{{\wt}a+\delta_{r}}\Res_{z_{0}}
\frac{(1+z_{0})^{{\wt}a-j-{r\over T}}}{z_{0}^{1+\delta_{r}}}Y(L(1)^{j}a,z_{0})
e^{L(1)}(-1)^{L(0)}b\nonumber\\
& &=\sum_{j=0}^{\infty}{1\over j!}(-1)^{{\wt}a+\delta_{r}}\Res_{z_{0}}
\frac{(1+z_{0})^{{\wt}(L(1)^{j}a)-1+{T-r\over
T}}}{z_{0}^{1+\delta_{r}}}Y(L(1)^{j}a,z_{0})
e^{L(1)}(-1)^{L(0)}b.\;\;\;\;\;\;
\end{eqnarray*}
By considering separately the cases $r=0$ and $r\ne 0$, we see that
the latter sum lies in $O_{g^{-1}}(V)$.  Thus $\phi$ maps $O_{g}(V)$ to
$O_{g^{-1}}(V)$. Since $\phi^{2}=1$ on $A(V^0)$, it is clear that
$\phi$ maps $O_{g}(V)$ onto $O_{g^{-1}}(V)$.  \qed

\section{Twisted modules}

We discuss the category of weak $g$-twisted $V$-modules (cf. [D1] and
[FFR]). As before $V$ is a
vertex operator algebra with automorphism $g$ of order $T$ and eigenspace
decomposition (\ref{g2.1}). We adopt standard notation in using
$W\{z\}$ to denote the
space of $W$-valued formal series in arbitrary real powers of $z$ for a vector
space $W.$
\begin{de} A {\em weak $g$-twisted $V$-module} $M$ is a vector space equipped
with a linear map
$$\begin{array}{l}
V\to (\End\,M)\{z\}\\
v\mapsto\displaystyle{ Y_M(v,z)=\sum_{n\in\Q}v_nz^{-n-1}\ \ \ (v_n\in
\End\,M)}
\end{array}$$
which satisfies the following for all $0\leq r\leq T-1,$ $u\in V^r$, $v\in V,$
$w\in M$,
\begin{eqnarray}
& &Y_M(u,z)=\sum_{n\in \frac{r}{T}+\Z}u_nz^{-n-1} \label{1/2}\\
& &v_lw=0\ \ \
\mbox{for}\ \ \ l>>0\label{vlw0}\\
& &Y_M({\bold 1},z)=1;\label{vacuum}
\end{eqnarray}
 \begin{equation}\label{jacobi}
\begin{array}{c}
\displaystyle{z^{-1}_0\delta\left(\frac{z_1-z_2}{z_0}\right)
Y_M(u,z_1)Y_M(v,z_2)-z^{-1}_0\delta\left(\frac{z_2-z_1}{-z_0}\right)
Y_M(v,z_2)Y_M(u,z_1)}\\
\displaystyle{=z_2^{-1}\left(\frac{z_1-z_0}{z_2}\right)^{-r/T}
\delta\left(\frac{z_1-z_0}{z_2}\right)
Y_M(Y(u,z_0)v,z_2)}.
\end{array}
\end{equation}
\end{de}

If $g=1$ this reduces to the definition of weak $V$-module as given in [DLM1].

One calls (\ref{jacobi}) the {\em twisted Jacobi identity,} and again it
reduces
to the ``untwisted'' Jacobi identity for $V$-modules if $g=1.$ Following
the arguments in the untwisted case (cf. [DL], [FHL], [FLM]) one can prove
that the twisted Jacobi identity is equivalent to the following
associativity and commutativity formulas:
\begin{eqnarray}\label{ea}
(z_{0}+z_{2})^{k+{r\over T}}Y_{M}(u,z_{0}+z_{2})Y_{M}(v,z_{2})w
=(z_{2}+z_{0})^{k+{r\over T}}Y_M(Y(u,z_0)v,z_2)w.
\end{eqnarray}
where $w\in M$ and $k$ is a nonnegative integer such that $z^{k+{r\over
T}}Y_{M}(u,z)w$ involves only positive powers of $z;$
\begin{eqnarray}
& &\ \ \ \  [Y_{M}(u,z_{1}),Y_{M}(v,z_{2})]\nonumber\\
& &=\Res_{z_{0}}z_2^{-1}\left(\frac{z_1-z_0}{z_2}\right)^{-r/T}
\delta\left(\frac{z_1-z_0}{z_2}\right)Y_M(Y(u,z_0)v,z_2).\label{ec}
\end{eqnarray}
Equating the coefficients of
$z_1^{-m-1}z_2^{-n-1}$ in (\ref{ec}) yields
\begin{eqnarray}\label{2.7}
[u_{m},v_{n}]=\sum_{i=0}^{\infty}
\left(\begin{array}{c}m\\i\end{array}\right)(u_{i}v)_{m+n-i}.
\end{eqnarray}

As in the untwisted case [DLM1], we may also deduce from
(\ref{1/2})-(\ref{jacobi}) the usual Virasoro algebra axioms, namely that
if $Y_M(\o,z)=\sum_{n\in\Z}L(n)z^{-n-2}$ then
\begin{equation}\label{g3.8}
[L(m),L(n)]=(m-n)L(m+n)+\frac{1}{12}(m^3-m)\delta_{m+n,0}(\mbox{rank}\,V)
\end{equation}
and
\begin{equation}\label{2.5}
\frac{d}{dz}Y_M(v,z)=Y_M(L(-1)v,z).
\end{equation}

The category of weak $g$-twisted $V$-modules has as its objects the weak
$g$-twisted $V$-modules and as morphisms those linear maps
$f: M\to W$ such that
\begin{equation}\label{g3.10}
fY_M(u,z)=Y_W(u,z)f
\end{equation}
for all $u\in V.$

\begin{de}\label{d3.2} A $g$-{\em twisted $V$-module} is
a weak $g$-twisted $V$-module $M$ which carries a
$\C$-grading induced by the spectrum of $L(0).$ That is, we have
\begin{equation}\label{g3.11}
M=\coprod_{\lambda \in{\C}}M_{\lambda}
\end{equation}
where $M_{\l}=\{w\in M|L(0)w=\l w\}.$ Moreover we require that
$\dim M_{\l}$ is finite and for fixed $\l,$ $M_{{n\over T}+\l}=0$
for all small enough integers $n.$
\end{de}

In this situation, if $w\in M_{\l}$ we refer to $\l$ as the {\em weight} of
$w$ and write $\l=\wt w.$ If $g=1$ then this defines a $V$-{\em
module} as used in [DLM1] and elsewhere. The totality of $g$-twisted
$V$-modules defines a full subcategory of the category of weak $g$-twisted
$V$-modules.

An important and related class of modules are the following.

\begin{de} An {\em admissible} $g$-twisted $V$-module
is a  weak $g$-twisted $V$-module $M$ which carries a
${1\over T}{\Z}_{+}$-grading
\begin{equation}\label{g3.12}
M=\oplus_{n\in\frac{1}{T}\Z}M(n)
\end{equation}
which satisfies the following
\begin{eqnarray}\label{g3.13}
v_{m}M(n)\subseteq M(n+\wt v-m-1)
\end{eqnarray}
for homogeneous $v\in V.$
\ed

The admissible $g$-twisted $V$-modules form a category with morphisms being
grade-preserving linear maps satisfying (\ref{g3.10}). Thus a {\em simple}
object in this category is an admissible $g$-twisted $V$-module $M$ such that
$0$ and $M$ are the only graded submodules.

We say that $V$ is $g$-{\em rational} if every admissible $g$-twisted
$V$-module is completely reducible, i.e., a direct sum of simple admissible
$g$-twisted modules.

If $g=1,$ these definitions reduce to the ``untwisted'' version used in [DLM1].
\begin{lem}\label{l3.4} There is a natural identification of the category
of $g$-twisted $V$-modules with a subcategory of the category
of admissible $g$-twisted
$V$-modules.
\el

\pf Let $M$ be a $g$-twisted $V$-module with decomposition into
$L(0)$-eigenspaces given by (\ref{g3.11}). For each $\l\in \C$ for which
$M_{\l}\ne 0,$ let $\l_0$ be the minimal element
of the set $\l+\frac{1}{T}\Z$ for
which $M_{\l_0}\ne 0.$ Note that $\l_0$ exists by definition \ref{d3.2}.
Let $\L$ be the set of all $\l_0$ so obtained, and for each
$n\in\frac{1}{T}\Z_+$ define
\begin{equation}\label{g3.14}
M(n)=\coprod_{\l\in\L}M_{n+\l}.
\end{equation}
It is clear that $M=\oplus_{n}M(n),$ while (\ref{g3.13}) follows from
the standard fact that $v_mM_{\l}\subset M_{\l+\wt v-m-1}.$

In this way we have identified $M$ as an admissible $g$-twisted $V$-module.
Moreover as a morphism $f$ in the category of $g$-twisted $V$-modules satisfies
(\ref{g3.10}) then it preserves $L(0)$-eigenspaces and hence also the
grading (\ref{g3.14}). The lemma follows. \qed

\begin{rem} We will establish later the less obvious fact that if $V$
is $g$-rational then the two categories of Lemma \ref{l3.4} share the
same simple objects.
\er
\begin{lem}\label{l3.7}
$M$ is  a simple weak $g$-twisted $V$-module. Then $M$ has  countable
dimension.
\end{lem}

\pf One knows (Proposition 2.4 of [DM2] or Lemma 6.1.1 of [L2]) that
$$M=\spa\{a_nu|a\in V,n\in\1t \Z_+\}$$
 for any non-zero $u\in M$.
Since $V$ has a countable basis, the lemma follows. \qed

This lemma is useful in the study of contragredient modules, which we now
discuss. If $M\!=\!\oplus_{n\in {1\over T}{\Z}_{+}}\!M(n)$
is an admissible $g$-twisted $V$-module, the contragredient module $M'$
is defined as follows:
\be{g3.15}
M'=\oplus_{n\in {1\over T}{\Z}_{+}}M(n)^{*}
\end{equation}
where $M(n)^*=\Hom_{\C}(M(n),\C).$ The vertex operator
$Y_{M'}(a,z)$ is defined for $a\in V$ via  by
\begin{eqnarray}\label{g3.16}
\langle Y_{M'}(a,z)f,u\rangle= \langle
f,Y_M(e^{zL(1)}(-z^{-2})^{L(0)}a,z^{-1})u\rangle
\end{eqnarray}

One can prove (cf. [FHL], [X]) the following:
\bl{l3.8} $(M',Y_{M'})$ is an admissible $g^{-1}$-twisted $V$-module.
\el

\section{The Lie algebra $V[g]$}

$V$\ \,continues\ \,to\ \,be\ a\ vertex\ operator algebra with automorphism $g$
of order
$T$ and eigenspace decomposition (\ref{g2.1}). Let $t$ be an indeterminate, and
consider the tensor product
\begin{equation}\label{g4.1}
{\cal L}(V)={\C}[t^{1\over T},t^{-{1\over T}}]\otimes V.
\end{equation}
Following [B] we give ${\C}[t^{1\over T},t^{-{1\over T}}]$ the structure
of vertex algebra with vertex operator
\begin{equation}\label{g4.2}
Y(f(t),z)g(t)=f(t+z)g(t)=\left(e^{z{d\over
dt}}f(t)\right)g(t).
\end{equation}
Then ${\cal L}(V)$ becomes a tensor product of vertex algebras and hence
itself a vertex algebra (cf. [DL], [FHL] and [L2]).

The action of $g$ naturally extends to that of a vertex
algebra automorphism
\be{g4.3}
g(t^{m}\otimes a)=\exp ({-2 \pi im\over
T})(t^{m}\otimes ga).
\end{equation}
Denote the space of $g$-invariants of this action by
${\cal L}(V,g);$ it is a vertex subalgebra of ${\cal L}(V).$ Of course
we have
\be{g4.4}
{\cal L}(V,g)=\oplus_{r=0}^{T-1}t^{r/T}\C[t,t^{-1}]\otimes V^r.
\end{equation}

The $L(-1)$ operator of ${\cal L}(V)$ and ${\cal L}(V,g)$ is given
by $D={d\over dt}\otimes 1+1\otimes L(-1),$ and as a consequence one
knows [B] that
\be{g4.5}
V[g]={\cal L}(V,g)/D{\cal L}(V,g)
\end{equation}
carries the structure of  Lie algebra with bracket
\be{g4.6}
[u+D{\cal L}(V,g),v+D{\cal L}(V,g)]=u_{0}v+D{\cal L}(V,g).
\end{equation}

As a matter of notation we use $a(q)$ to denote the image of $t^q\otimes a\in
{\cal L}(V,g)$ in $V[g].$ An easy computation from the definitions yields
\bl{l4.1} Let $a\in V^r,$ $v\in V^s$ and $m,n\in\Z.$ Then

(i) $[\omega(0),a(m+{r\over T})]=-\left(m+{r\over T}\right)a(m-1+{r\over T}),$

(ii) $[a(m+{r\over T}), b(n+{s\over T})]=
\sum_{i=0}^{\infty}{m+{r\over T}\choose i}a_ib(m+n+{r+s\over T}-i),$

(iii) $\o(0)$ and ${\bold 1}(-1)$ both lie in the center of $V[g].$
\el

We can introduce a ${1\over T}{\Z}$-gradation on ${\cal L}(V)$ by defining,
for homogeneous $a\in V,$
\be{g4.7}
\deg(t^n\otimes a)=\wt a-n-1.
\end{equation}
As $D$ increases degree by 1 then $D{\cal L}(V,g)$ is a graded
subspace of ${\cal L}(V,g),$ so that there is a naturally induced
$\frac{1}{T}\Z$-gradation on $V[g].$ Let $V[g]_n$ denote the degree $n$
subspace of $V[g].$ After Lemma \ref{l4.1} $V[g]$ is in fact
a $\frac{1}{T}\Z$-graded Lie algebra and we have a triangular decomposition
\be{g4.8}
V[g]=V[g]_{+}\oplus V[g]_0\oplus V[g]_{-}
\end{equation}
where we have set
$\displaystyle{V[g]_{\pm}=\sum_{0<n\in {1\over T}{\Z}}V[g]_{\pm n}}$.

Note that $V[g]_0$ is spanned by elements of the form
$a(\wt a-1)$ for homogeneous $a\in V^0.$ The bracket is given
by
\begin{eqnarray}\label{ec0}
[a({\wt}a-1),b({\wt}b-1)]=\sum_{j=0}^{\infty}{{\wt}a-1\choose
j}a_jb({\wt}(a_{j}b)-1).
\end{eqnarray}
as we see from Lemma \ref{l4.1} (ii).

Following [Z] we set $o(a)=a(\wt a-1)$ for homogeneous $a\in V^0.$ So we
have a linear map
\begin{eqnarray}
& & V^0\to V[g]_0,\nonumber\\
& & a\mapsto o(a)\ (a\ {\rm homogeneous}).\label{g4.10}
\end{eqnarray}
The kernel of the map is precisely $(L(-1)+L(0))V^0,$ and (\ref{g4.10})
induces an isomorphism of Lie algebras $V^0/(L(-1)+L(0))V^0\cong V[g]_0$
where the bracket on the quotient of $V^0$ is as described via
$$[a,b]=\sum_{j\geq 0}{{\wt}a-1\choose j}a_{j}b.$$

\bl{l4.2} Let $A_g(V)_{Lie}$ be the Lie algebra of the associative algebra
$A_g(V)$ (cf. Section 2), so that $[u,v]=u*_gv-v*_gu.$ Then the map
$o(a)\mapsto a+O_g(V)$ is a Lie algebra epimorphism $V[g]_0\to A_g(V)_{Lie}.$
\el

\pf Let $I=O_g(V)\cap V^0.$ We have from Lemma 2.1.1 of [Z] that
$(L(-1)+L(0))V^0\subset O(V^0)\subset I,$ so we have surjective linear maps
\be{g4.11}
V[g]_0\cong V^0/(L(-1)+L(0))V^0\to A(V^0)\to V^0/I\cong A_g(V).
\end{equation}
Use Lemma \ref{l2.2} (ii) to see that if $a,b\in V^0$ are homogeneous then
$$a*b-b*a\equiv \Res_{z}(1+z)^{\wt a-1}Y(a,z)b=\sum_{i=0}^{\infty}
{\wt a-1\choose i}a_ib\ (\m O(V^0)).$$
In view of the results following (\ref{g4.10}) this says that the map
from $V^0/(L(-1)+L(0))V^0$ to $A(V^0)$ in (\ref{g4.11}) is a Lie algebra
morphism, and the lemma follows. \qed

\br{r4.3}Observe that $V[g]$ contains $V^0[1]=\C[t,t^{-1}]\otimes
V^0/D(\C[t,t^{-1}]\otimes V^0)$ as a graded Lie subalgebra.
\er

\section{The functor $\O$}

In this section we construct a  functor $\Omega$
from the category of admissible $g$-twisted $V$-modules to the category
of $A_g(V)$-modules. We retain previous notation.

The connection between the Lie algebra $V[g]$ and weak $g$-twisted modules is
the following:
\bl{l5.1} Let $M$ be a weak $g$-twisted $V$-module. The map $a(m)\mapsto a_m$
defines a representation of the Lie algebra $V[g]$ on $M.$
\el

\pf The bracket of elements $a(m),b(n)$ in $V[g]$ $(m,n\in\frac{1}{T}\Z)$ is
given in Lemma \ref{l4.1} (ii), and that of $a_m,b_n$ is given in
(\ref{2.7}).
Comparing, we see that it suffices to show that the map $a(m)\mapsto
a_m$ is well-defined. Let $t^m\otimes a\in{\cal L}(V,g).$ Then
$D(t^m\otimes a)=mt^{m-1}\otimes a+t^{m}\otimes L(-1)a \mapsto
ma_{m-1}+(L(-1)a)_m=0,$ the latter equality by (\ref{2.5}).
\qed

\bl{l5.2} Let $M$ be a weak $g$-twisted $V$-module which carries a
$\frac{1}{T}\Z_+$-grading. Then $M$ is an admissible $g$-twisted
$V$-module if, and only if, $M$ is a $\frac{1}{T}\Z_+$-graded module for
the grade Lie algebra $V[g].$
\el

\pf Let $a\in V$ be homogeneous and $a(m)\in V[g],$with
$M(n)$ the $n$-th graded piece of $M$ (cf. (\ref{g3.12}). The condition
that $M$ is graded module for $V[g]$ is this: $a(m)M(n)\subset
M(n+\deg a(m))=M(n+\wt a-m-1).$ Using the representation described in Lemma
\ref{l5.1}, this is precisely the condition (ii) of Lemma \ref{l4.1}
 required to make $M$
an admissible $g$-twisted module. \qed

Recalling the decomposition (\ref{g4.8}) of $V[g],$ consider
a module $W$ for the Lie algebra $V[g].$ We let $\O(W)$ denote the space
of ``lowest
weight vectors,'' that is
\begin{eqnarray}\label{g5.1}
\Omega (W)=\{u\in W|V[g]_{-}u=0\}.
\end{eqnarray}
Similarly, for a $V^0[1]$-module $W_0$ (cf. Remark \ref{r4.3})
we set
\begin{equation}\label{g5.2}
\Omega^0(W_0)=\{u\in W_0|V^0[1]_{-}u=0\}.
\end{equation}
Then $\O(W)$ and $\O^0(W_0)$ are modules for the Lie algebras $V[g]_0$ and
$V^0[1]_0$ respectively. Furthermore it is obvious from Remark \ref{r4.3} that
we have $\O(W)\subset \O^0(W)$ for a $V[g]$-module $W.$

\bt{t5.3} Suppose that $M$ is  weak $g$-twisted $V$-module. Then there is
a representation of the associative algebra $A_g(V)$ on $\O(M)$ induced by
the map $a\mapsto o(a)$ for homogeneous $a\in V^0$ (cf. (\ref{g4.10}).
\et

\pf We start by remarking that if $g=1$
then this result has been established by Zhu [Z]. Although he works in a less
general situation, one easily verifies that his proof goes through in the
present situation.

We make use of this as follows: the theorem is correct as applied to the
action of $A(V^0)$ on $\O^0(M).$ So we are reduced to proving that
$\O(M)$ is an $A(V^0)$-submodule of $\O^0(M)$ on which $O_g(V)\cap V^0$
acts trivially.

To show that $\O(M)$ is $A(V^0)$-stable, pick $u\in \Omega(M)$ and
homogeneous $a\in V^0.$ We must show that
$o(a)u$ is annihilated by all $b(n)\in V[g]$ where  $b\in V$ is homogeneous
and $\deg b(n)<0,$  that is, $b_no(a)u=0.$

We have, using (\ref{2.7}),
$$ b_no(a)u=b_na_{\wt a-1}u=a_{\wt a-1}b_nu+\sum_{i=0}^\infty{n\choose i}
(b_ia)_{m+n-i-1}u.$$
Now $b_nu=0$ since $b(n)\in V[g]_-$ annihilates $\O(M).$ And since
$b_ia(n+\wt a -i-1)$ has degree equal to $(\wt b+\wt a-i-2)-(n+\wt a -i-1)-1=
\wt b-n-1=\deg b(n)<0$ then each term $(b_ia)_{\wt a+n-i-1}u=0.$ So indeed
$\O(M)$ is an $A(V^0)$-module.

It remains to prove that for any $a\in O_g(V)\cap V^0$,
$o(a)$ acts as zero on $\Omega(M)$. If $a\in O(V^0)$ then
$o(a)=0$ since $\Omega(M)$ is an $A(V^0)$-module.
Suppose, then, that
$$a=\Res_{z}\frac{(1+z)^{{\wt}c-1+{r\over T}}}{z}Y(u,z)v$$
with $u\in V^r, v\in V^{T-r}, 1\le r\le T-1$. For any $w\in \Omega(M)$, using
a property of the
delta-function we can rewrite the Jacobi identity (\ref{jacobi}) as follows:
\begin{equation}
\begin{array}{c}\label{e2.23}
\displaystyle{z^{-1}_1\delta\left(\frac{z_0+z_2}{z_1}\right)
Y_M(u,z_1)Y_M(v,z_2)w-z^{-1}_0\delta\left(\frac{z_2-z_1}{-z_0}\right)
Y_M(v,z_2)Y_M(u,z_1)w}\\
\displaystyle{=z_1^{-1}\left(\frac{z_2+z_0}{z_1}\right)^{r/T}
\delta\left(\frac{z_2+z_0}{z_1}\right)
Y_M(Y(u,z_0)v,z_2)w.}
\end{array}
\end{equation}
Apply $\Res_{z_{1}}z_{1}^{{\wt}u-1+{r\over T}}$ to
(\ref{e2.23}) to obtain the following
twisted associativity:
\begin{eqnarray}\label{e2.24}
& &\ \ \ \Res_{z_{1}}z_{1}^{{\wt}u-1+{r\over
T}}z^{-1}_1\delta\left(\frac{z_0+z_2}{z_1}\right)
Y_M(u,z_1)Y_M(v,z_2)w\nonumber\\
& &=(z_{2}+z_{0})^{{\wt}u-1+{r\over T}}Y_{M}(Y(u,z_{0})v,z_{2})w.
\end{eqnarray}
Take $\Res_{z_{0}}\Res_{z_{2}}z_{0}^{-1}z_{2}^{{\wt}v-{r\over T}}$ of
(\ref{e2.24}) to get
\begin{eqnarray}
& &0=\Res_{z_{0}}\Res_{z_{2}}z_{0}^{-1}z_{2}^{{\wt}v-{r\over T}}
(z_{2}+z_{0})^{{\wt}u-1+{r\over T}}Y_{M}(Y(u,z_{0})v,z_{2})w\nonumber\\
& &\ \ \ \ =\sum_{i=0}^{\infty}\left(\begin{array}{c}{\wt}u-1+{r\over
T}\\i\end{array}\right)
\Res_{z_{2}}z_{2}^{{\wt}u+{\wt}v-i-1}Y_{M}(u_{i-1}v,z_{2})w\nonumber\\
& &\ \ \ \ =\sum_{i=0}^{\infty}\left(\begin{array}{c}{\wt}u-1+{r\over
T}\\i\end{array}\right)o(u_{i-1}v)w
\nonumber\\
& &\ \ \ \ =o\left(\Res_{z}\frac{(1+z)^{{\wt}u-1+{r\over
T}}}{z}Y_{M}(u,z)v\right)w\nonumber\\
& &\ \ \ \ =o(a)w
\end{eqnarray}
as required. \qed

Because our constructions are natural, it is evident that $\Omega$ is
a covariant functor from the category of
weak $g$-twisted $V$-modules
to the category of $A_{g}(V)$-modules. To be more precise, if $f:M\to N$
is a morphism in the first category (cf. (\ref{g3.10}) we define $\O(f)$
to be the restriction of $f$ to $\O(M)$. With an obvious notation,
(\ref{g3.10}) says that $fa_m^M=a_m^Nf$ for $a\in V$ and
$m\in\frac{1}{T}\Z.$ Then $f$ induces a morphism of $V[g]$-modules $M\to N$
by Lemma \ref{l5.1}. Moreover
$\O(f):\O(M)\to\O(N).$ Now Theorem \ref{t5.3} implies that $\O(f)$ is a
morphism of $A_g(V)$-modules.

We turn to a consideration of admissible $g$-twisted $V$-modules in this
context. Let $M$ be such a module. As long as $M\ne 0,$ then some
$M(n)\ne 0$ (cf. (\ref{g3.12}), and it is no loss to shift the grading
so that in fact $M(0)\ne 0.$ If $M=0,$ let $M(0)=0.$ With these
conventions we prove
\bp{l2.9} Suppose that $M$ is a simple  admissible $g$-twisted $V$-module.
Then the following hold

(i) $\Omega(M)=M(0).$

(ii) $\O(M)$ is a simple $A_{g}(V)$-module.
\ep

\pf Note that Lemma \ref{l5.2} is available in this situation. An
easy argument shows that $\O(M)$ is a graded subspace of $M.$ That is
\begin{equation}\label{g5.6}
\O(M)=\oplus_{n\in\frac{1}{T}\Z}\O(M)\cap M(n).
\end{equation}
Set $\O(n)=\O(M)\cap M(n).$ It is clear that $M(0)\subset \O(M).$
In order to prove (i) we must show that $\O(n)=0$ if $n>0.$
We use the PBW theorem to do this.

Let $U(\cdot)$ denote universal enveloping algebra. If $\O(n)\ne 0$ then
because $M$ is simple we have
\be{g5.7}
M=U(V[g])\Omega(n)=U(V[g]_{+})\Omega(n),
\end{equation}
the latter equality thanks to the triangular decomposition of $V[g]$
(\ref{g4.8}). Equation (\ref{g5.7}) tells us that the lowest
degree of $M$ is no less than $n,$ so we must have $n=0$ by our convention.

To prove (ii) let $U$ be
any nonzero $A_{g}(V)$-submodule of $M(0)$. Then $U$ is annihilated
by $V[g]_-$ and stable under $V[g]_0$ (Lemma \ref{l4.2} and Theorem
\ref{t5.3}). So again the PBW theorem yields
\be{g5.8}
M=U(V[g])U=U(V[g]_{+})U=U\oplus U(V[g]_{+})V[g]_{+}U.
\end{equation}
This implies that $U=M(0),$ and (ii) follows. \qed

%\bl{l5.5} Let $M$ be a weak $g$-twisted $V$-module. Then $U(V[g])\O(M)$
%is an admissible $g$-twisted submodule of $M.$
%\el

%\pf We see from Lemma \ref{l5.1} that the $V[g]$-submodules of $M$ are
%precisely the $g$-twisted $V$-submodules of $M.$ So $U(V[g])\O(M)$ is
%certainly a $g$-twisted $V$-submodule of $M;$ it is admissible by
%Lemma \ref{l5.2} if we endow $\O(M)$ with degree 0. Remark that if $M$
%is itself admissible then $\O(M)$ is a $\frac{1}{T}\Z_+$-grade
%subspace, in which case $U(V[g])\O(M)$ is an admissible submodule. \qed

\section{Generalized Verma modules and the functor $L$}

We consider the possibility of constructing admissible $g$-twisted $V$-modules
from a given $A_g(V)$-module $U,$ say. We show that there is a universal
way to do this. Moreover a certain quotient of the universal object is an
admissible $g$-twisted $V$-module $L(U)$ and $L$ defines a functor which
is right inverse to the functor $\O.$ Notation is as before.

Given the $A_{g}(V)$-module $U$, it is a fortiori a module for $A_g(V)_{Lie}.$
Thanks to Lemma \ref{l4.2} we can lift $U$ to a module for the Lie algebra
$V[g]_0,$ and then to one for $V[g]_{-}\oplus V[g]_{0}$ by letting $V[g]_-$
act trivially. Set $P=V[g]_{-}\oplus V[g]_{0}$ and define
\be{g6.1}
M(U)=\Ind_{P}^{V[g]}(U)=U(V[g])\otimes_{U(P)} U.
\end{equation}
If we give $U$ degree 0, the $\frac{1}{T}\Z$-gradation of $V[g]$ lifts to
$M(U)$ which thus becomes a  $\frac{1}{T}\Z_+$-graded module for $V[g].$
The PBW theorem implies that $M(U)(n)=U(V[g]_+)_nU$ and in
particular $M(U)(0)=U.$

Taking our cue from Lemma \ref{l5.1}, we define for $v\in V^r,$
\be{g6.2}
Y_{M(U)}(v,z)=\sum_{n\in\frac{1}{T}\Z}v(m)z^{-1-m}
\end{equation}
where we convene that $v(m)=0$ unless $v(m)=t^m\otimes v$ lies in $V[g].$
Then $Y_{M(U)}(v,z)$ satisfies condition (\ref{1/2}). Moreover (\ref{vlw0})
and (\ref{vacuum}) are easily confirmed.

Next, thanks to Lemma \ref{l4.1} (ii), we see that the identity (\ref{2.7})
 holds. Thus in order to establish the $g$-twisted Jacobi identity for the
action (\ref{g6.2}) on $M(U)$ it would be enough to also establish (\ref{ea}).
In general, however, this is false. Instead we have to divide out by
the desired relations.

Precisely, let $W$ be the subspace of $M(U)$ spanned linearly by the
coefficients of
\begin{eqnarray}\label{g6.3}
(z_{0}+z_{2})^{{\wt}a-1+\delta_r+{r\over
T}}Y(a,z_{0}+z_{2})Y(b,z_{2})u-(z_{2}+z_{0})^{{\wt}a-1+\delta_r+{r\over T}}
Y(Y(a,z_{0})b,z_{2})u
\end{eqnarray}
for any homogeneous $a\in V^r,b\in V,$ $u\in U$.
We set
\be{g6.4}
\bar M(U)=M(U)/U(V[g])W.
\end{equation}

In order to prove the first main result of this section we need the
following proposition.
\bp{p6.1}  Let $M$ be a $V[g]$-module such that there is a
 subspace $U$ of $M$ satisfying the
following conditions:

(i) $M=U(V[g])U;$

(ii) For any $a\in V^r$ and $u\in U$ there is a positive
integer $k$ such that
\begin{eqnarray}\label{d1}
(z_{0}+z_{2})^{k+{r\over T}}Y(a,z_{0}+z_{2})Y(b,z_{2})u=
(z_{0}+z_{2})^{k+{r\over T}}Y(Y(a,z_{0})b,z_{2})u
\end{eqnarray}
for any $b\in V$. Then $M$ is a weak $V$-module.
\ep

\pf. It suffices to show that (\ref{d1}) holds for all $u\in M.$
Let $X$ consist of those $u\in M$ such that for any $a\in
V^r$ the associativity (\ref{d1}) holds
for any $b\in V$. We must show that $X=M$. Since $X$
contains $U$ by (ii) and $M$ is generated by $U$ it is equivalent to
show that $V[g]X\subset X.$

Let $u\in X, c\in V$ and $n\in\1t\Z.$
Let $k_{1}$ be a positive integer such that $c_{i}a=0$
for $i\ge
k_{1}$. Since $u\in X$, there is a positive integer $k_{2}$ such that
\begin{eqnarray}
& &(z_{0}+z_{2})^{k_{2}+{r+s\over T}}Y(c_{i}a,z_{0}+z_{2})Y(b,z_{2})u
=(z_{2}+z_{0})^{k_{2}+{r+s\over T}}Y(Y(c_{i}a,z_{0})b,z_{2})u,\label{3.18}\\
& &(z_{0}+z_{2})^{k_{2}+{r+s\over T}}Y(a,z_{0}+z_{2})Y(c_{i}b,z_{2})u
=(z_{2}+z_{0})^{k_{2}+{r+s\over
T}}Y(Y(a,z_{0})c_{i}b,z_{2})u\hspace{1cm}\label{3.19}
\end{eqnarray}
for any nonnegative integer $i$. Let $k$ be a positive integer such that
$k+{r\over T}+n-k_{1}>k_{2}+{r+s\over T}$. Use (\ref{3.18}) and
(\ref{3.19}) and the equality (a consequence of (\ref{2.7}))
\begin{eqnarray}\label{3.6'}
[a_m,Y(b,z_{2})]
=\sum_{i=0}^{\infty}{m\choose i}z_2^{m-i}Y(a_ib,z_2)
\end{eqnarray}
to obtain
\begin{eqnarray*}
& &\ \ \ \ \ (z_{0}+z_{2})^{k+{r\over T}}Y(a,z_{0}+z_{2})Y(b,z_{2})c_{n}u\\
& &=(z_{0}+z_{2})^{k+{r\over T}}c_{n}Y(a,z_{0}+z_{2})Y(b,z_{2})u\\
& &\ \ \ \ -\sum_{i=0}^{\infty}{n\choose i}
(z_{0}+z_{2})^{k+{r\over T}+n-i}Y(c_{i}a,z_{0}+z_{2})Y(b,z_{2})u\\
& &\ \ \ \  -\sum_{i=0}^{\infty}{n\choose i}
z_{2}^{n-i}(z_{0}+z_{2})^{k+{r\over T}}Y(a,z_{0}+z_{2})Y(c_{i}b,z_{2})u\\
& &=(z_{0}+z_{2})^{k+{r\over T}}c_{n}Y(Ya,z_{0})b,z_{2})u\\
& &\ \ \ \ -\sum_{i=0}^{\infty}{n \choose i}
(z_{2}+z_{0})^{k+{r\over T}+n-i}Y(Y(c_{i}a,z_{0})b,z_{2})u\\
& &\ \ \ \ -\sum_{i=0}^{\infty}{n \choose i}
z_{2}^{n-i}(z_{2}+z_{0})^{k+{r\over T}}Y(Y(a,z_{0})c_{i}b,z_{2})u\\
& &=(z_{0}+z_{2})^{k+{r\over T}}c_{n}Y(Ya,z_{0})b,z_{2})u\\
& &\ \ \ \ -\sum_{i=0}^{\infty}{n \choose i}
(z_{2}+z_{0})^{k+{r\over T}+n-i}Y(Y(c_{i}a,z_{0})b,z_{2})u\\
& &\ \ \ \ -\sum_{i=0}^{\infty}{n \choose i}
z_{2}^{n-i}(z_{2}+z_{0})^{k+{r\over T}}Y(c_{i}Y(a,z_{0})b,z_{2})u\\
& &\ \ \ \ +\sum_{i=0}^{\infty}\sum_{j=0}^{\infty}{n \choose j}{j \choose i}
z_{2}^{n-i}(z_{2}+z_{0})^{k+{r\over T}}z_{0}^{i-j}Y(Y(c_{j}a,z_{0})b,z_{2})u\\
& &=(z_{0}+z_{2})^{k+{r\over T}}c_{n}Y(Ya,z_{0})b,z_{2})u\\
& &\ \ \ \ -\sum_{i=0}^{\infty}{n\choose i}
(z_{2}+z_{0})^{k+{r\over T}+n-i}Y(Y(c_{i}a,z_{0})b,z_{2})u\\
& &\ \ \ \ -\sum_{i=0}^{\infty}{n\choose i}
z_{2}^{n-i}(z_{2}+z_{0})^{k+{r\over T}}Y(c_{i}Y(a,z_{0})b,z_{2})u\\
& &\ \ \ \ +\sum_{j=0}^{\infty}\sum_{i=j}^{\infty}{n\choose j}{n-j\choose i-j}
z_{2}^{n-i}(z_{2}+z_{0})^{k+{r\over T}}z_{0}^{i-j}
Y(Y(c_{j}a,z_{0})b,z_{2})u \nonumber\\
& &=(z_{0}+z_{2})^{k+{r\over T}}c_{n}Y(Ya,z_{0})b,z_{2})u\\
& &\ \ \ \ -\sum_{i=0}^{\infty}{n\choose i}
z_{2}^{n-i}(z_{2}+z_{0})^{k+{r\over T}}Y(c_{i}Y(a,z_{0})b,z_{2})u\\
& &=(z_{0}+z_{2})^{k+{r\over T}}c_{n}Y(Ya,z_{0})b,z_{2})u\\
& &\ \ \  \ -(z_{2}+z_{0})^{k+{r\over T}}[c_n, Y(Y(a,z_{0})b,z_{2})]u\\
& &=(z_{2}+z_{0})^{k+{r\over T}}Y(Y(a,z_{0})b,z_{2})c_{n}u,
\end{eqnarray*}
as desired. \qed

The first main result of this section is the following:
\bt{t6.1} $\bar M(U)$ is an admissible $g$-twisted $V$-module with
$\bar M(U)(0)=U$ and with the
following universal property: for any weak $g$-twisted $V$-module $M$
and any $A_g(V)$-morphism $\phi: U\to \O(M),$ there is a unique morphism
$\bar\phi: \bar M(U)\to M$ of weak $g$-twisted $V$-modules which
extends $\phi.$
\et

\pf Clearly $\bar M(U)$ satisfies the conditions placed on $M$ in Proposition
\ref{p6.1}. We see that $\bar M(U)$ is a weak $g$-twisted $V$-module and
therefore an admissible $g$-twisted $V$-module. The universal property of
$\bar M(U)$ follows from its construction. A proof of that $\bar M(U)(0)=U$
will be given after Proposition \ref{p3.5}. \qed

We can now state the second main result of this section
\bt{t6.3} $M(U)$ has a unique maximal graded $V[g]$-submodule $J$
with the property that $J\cap U=0.$ Then $L(U)=M(U)/J$ is an admissible
$g$-twisted $V$-module satisfying $\O(L(U))\cong U.$

$L$ defines a functor from the category of
$A_{g}(V)$-modules to
the category of admissible $g$-twisted $V$-modules such that $\O\circ L$
is naturally equivalent to the identity.
\et

In the following we let $U^*=\Hom_{\C}(U,\C)$ and extend
$U^*$ to $M(U)$ by letting $U^*$ annihilate $\oplus_{n>0}M(U)(n).$
Then one easily shows
\bl{l6.4} Let $J$ be the graded $V[g]$-submodule of $M(U)$ maximal subject
to $J\cap U=0.$ Then
$$ J=\{v\in M(U)|\langle u',xv\rangle=0\ {\rm for\ all}\ u'\in
U^{*},\ {\rm all}\ x\in U(V[g])\}.$$
\el

The main point in the proof of the Theorems is to show that
$U(V[g])W\subset J.$ The next three results are devoted to this goal.

\begin{prop}\label{p3.3} The following hold
for all homogeneous $a\in V^r, b\in V,$ $u'\in U^{*}, u\in U,j\in {\Z}_{+},$
\begin{eqnarray}\label{ea3}
& &\ \ \ \langle u',(z_{0}+z_{2})^{{\wt}a-1+\delta_{r}+{r\over
T}+j}Y_{M(U)}(a,z_{0}+z_{2})Y_{M(U)}(b,z_{2})u\rangle
\nonumber\\
& &=\langle u',(z_{2}+z_{0})^{{\wt}a-1+\delta_{r}+{r\over
T}+j}Y_{M(U)}(Y(a,z_{0})b,z_{2})u\rangle.
\end{eqnarray}
\ep
\br{r6.6} In the following we habitually drop subscripts attached to $Y,$ which
should cause no confusion.
\er

{\bf Proof of Proposition \ref{p3.3}:}
By linearity we may take $b\in V^s$ homogeneous. Suppose that $v(m)\in V[g]$
and $u\in U.$  Notice that $Y(a,z)=\sum_{n\in \frac{r}{T}+\Z}a(n)z^{-n-1}.$
It follows that if $r+s\not\equiv 0$
(\m $T$) then all coefficients $z_0^iz_2^j$
of $Y(a,z_{0}+z_{2})Y(b,z_{2})u$ and $Y(Y(a,z_0)b,z_2)$ are either
0 or lie in $M(U)(n)$ for $n>0.$
So in this case both sides of (\ref{ea3}) are zero since $u'$
annihilates $M(U)(n)$ for $n>0.$

So we may assume that $r+s\equiv 0$ (\m $T$). If $r=s=0$ we may appeal either
to [Z] or [L2] as this is essentially the case $g=1.$ So from now on we
take $1\le r\le T-1,$ $s=T-r.$ Note that $\delta_r=0$ in this case. We proceed
in several lemmas.
\bl{l6.7} For any $i,j\in \Z_{+},$
\begin{eqnarray*}
& &\ \ \Res_{z_{0}}z_{0}^{-1+i}(z_{0}+z_{2})^{{\wt}a-1+{r\over T}+j}
\langle u',Y(a,z_{0}+z_{2})Y(b,z_{2})u\rangle\nonumber\\
& &=\Res_{z_{0}}z_{0}^{-1+i}(z_{2}+z_{0})^{{\wt}a-1+{r\over T}+j}
\langle u',Y(Y(a,z_{0})b,z_{2})u\rangle.
\end{eqnarray*}
\el

\pf Since $j\geq 0$ then $a(\wt a-1+\frac{r}{T}+j)$ lies in $V[g]_-$ and hence
annihilates $u.$ Then for all $i\in\Z_+$ we get
\begin{eqnarray}\label{3.5}
\Res_{z_{1}}(z_{1}-z_{2})^{i}z_{1}^{{\wt}a-1+{r\over T}+j}
Y(b,z_{2})Y(a,z_{1})u=0.
\end{eqnarray}
Note that Lemma \ref{l4.2} (ii) is equivalent to
\begin{eqnarray}\label{3.6''}
[Y(a,z_{1}),Y(b,z_{2})]
=\Res_{z_{0}}z_2^{-1}\left(\frac{z_1-z_0}{z_2}\right)^{-r/T}
\delta\left(\frac{z_1-z_0}{z_2}\right)(Y(u,z_0)v,z_2)
\end{eqnarray}
(cf. (\ref{ec}).

Using (\ref{3.5}) and (\ref{3.6''})
we obtain:
\begin{eqnarray}
& &\ \ \ \Res_{z_{0}}z_{0}^{i}(z_{0}+z_{2})^{{\wt}a-1+j+{r\over T}}
Y(a,z_{0}+z_{2})Y(b,z_{2})u\nonumber\\
& &= \Res_{z_{1}}(z_{1}-z_{2})^{i}z_{1}^{{\wt}a-1+{r\over T}+j}
Y(a,z_{1})Y(b,z_{2})u\nonumber\\
& &= \Res_{z_{1}}(z_{1}-z_{2})^{i}z_{1}^{{\wt}a-1+{r\over T}+j}
Y(a,z_{1})Y(b,z_{2})u\nonumber\\
& &\ \ \ -\Res_{z_{1}}(z_{1}-z_{2})^{i}z_{1}^{{\wt}a-1+{r\over T}+j}
Y(b,z_{2})Y(a,z_{1})u\nonumber\\
& &= \Res_{z_{1}}(z_{1}-z_{2})^{i}z_{1}^{{\wt}a-1+{r\over T}+j}
[Y(a,z_{1}),Y(b,z_{2})]u\nonumber\\
& &= \Res_{z_{0}}\Res_{z_{1}}(z_{1}-z_{2})^{i}z_{1}^{{\wt}a-1+{r\over T}+j}
z_{2}^{-1}\delta\left(\frac{z_{1}-z_{0}}{z_{2}}\right)
\left(\frac{z_{1}-z_{0}}{z_{2}}\right)^{-{r\over T}}
Y(Y(a,z_{0})b,z_{2})u\nonumber\\
& &=\Res_{z_{0}}\Res_{z_{1}}z_{0}^{i}z_{1}^{{\wt}a-1+{r\over T}+j}
z_{1}^{-1}\delta\left(\frac{z_{2}+z_{0}}{z_{1}}\right)
\left(\frac{z_{2}+z_{0}}{z_{1}}\right)^{{r\over T}}Y(Y(a,z_{0})b,z_{2})u
\nonumber\\
& &=\Res_{z_{0}}z_{0}^{i}(z_{2}+z_{0})^{{\wt}a-1+{r\over T}+j}
Y(Y(a,z_{0})b,z_{2})u.\label{3.6}
\end{eqnarray}
Thus lemma \ref{l6.7} holds if $i\geq 1,$ and we may now assume
$i=0.$

Next use (\ref{3.6}) to calculate that
\begin{eqnarray}
& &\Res_{z_{0}}z_{0}^{-1}(z_{0}+z_{2})^{{\wt}a-1+{r\over T}+j}
\langle u',Y(a,z_{0}+z_{2})Y(b,z_{2})u\rangle\nonumber\\
&=&\sum_{k=0}^{\infty}\left(\begin{array}{c}j\\k\end{array}\right)
\Res_{z_{0}}z_{0}^{k-1}z_{2}^{j-k}(z_{0}+z_{2})^{{\wt}a-1+{r\over T}}
\langle u',Y(a,z_{0}+z_{2})Y(b,z_{2})u\rangle\nonumber\\
&=&\sum_{k=1}^{\infty}\left(\begin{array}{c}j\\k\end{array}\right)
\Res_{z_{0}}z_{0}^{k-1}z_{2}^{j-k}(z_{2}+z_{0})^{{\wt}a-1+{r\over T}}
\langle u',Y(Y(a,z_{0})b,z_{2})u\rangle\nonumber\\
& & \ \ \ + \Res_{z_{0}}z_{0}^{-1}z_2^j(z_{2}+z_{0})^{{\wt}a-1+{r\over T}}
\langle u',Y(a,z_0+z_2)Y(b,z_{2})u\rangle.\label{g6.9}
\end{eqnarray}

We claim that we have
\begin{eqnarray}\label{g6.10}
& &\Res_{z_{0}}z_{0}^{-1}(z_{0}+z_{2})^{{\wt}a-1+{r\over T}}
\langle u',Y(a,z_{0}+z_{2})Y(b,z_{2})u\rangle =0
\end{eqnarray}
and
\begin{eqnarray}\label{g6.11}
\Res_{z_{0}}z_{0}^{-1}(z_{2}+z_{0})^{{\wt}a-1+{r\over T}}
\langle u',Y(Y(a,z_{0})b,z_{2})u\rangle=0.
\end{eqnarray}
If so, then Lemma \ref{l6.7} follows from (\ref{g6.9}).

To see (\ref{g6.10}), note that $b_{\wt b-1-n}\in V[g]_-$ if $n<0,$ so
that $b_{\wt b-1-n}u=0$ in this case. Then we see that
\begin{eqnarray}
& &\ \ \ \ \Res_{z_{0}}z_{0}^{-1}(z_{0}+z_{2})^{{\wt}a-1+{r\over T}}
\langle u',Y(a,z_{0}+z_{2})Y(b,z_{2})u\rangle\nonumber\\
& &=\<u',\sum_{k\in\Z_+}a(\wt a-2-k+\frac{r}{T})\sum_{n\geq 0}b(\wt
b-1-n)z_2^{-\wt b+n+k}u\>.\label{g6.12}
\end{eqnarray}
Since all operators of the form $a(q)$ occurring in (\ref{g6.12})
lie in $V[g]_+$ then all components lie in $\oplus_{n>0}M(U)(n)$
and hence are annihilated by $u'.$ So (\ref{g6.10}) follows.

Completely analogous calculation shows that
$$\Res_{z_{0}}z_{0}^{-1}(z_{2}+z_{0})^{{\wt}a-1+{r\over T}}
\langle u',Y(Y(a,z_{0})b,z_{2})u\rangle$$
 is a sum of terms with components
annihilated by $u',$ plus the expression
\be{g6.13}
z_{2}^{-{\wt }b-1+{r\over T}}\sum_{k\in\Z_+}{\wt a-1+{r\over T}\choose k}
(a_{k-1}b)(\wt(a_{k-1}b)-1)u
\end{equation}
In the notation of (\ref{g4.10}), (\ref{g6.13}) is equal to
$$z_{2}^{-{\wt }b-1+{r\over T}}\sum_{k\in\Z_+}{\wt a-1+{r\over T}\choose k}
o(a_{k-1}b)u$$
(note that $a_{k-1}b\in V^0$), and from (\ref{g2.2}) and Lemma \ref{l4.2}
the sum is precisely the action of $a\circ_gb$ on $u.$ Since
$a\circ_gb\in O_g(V)$ we have $(a\circ_gb)u=0$ (cf. Theorem \ref{t5.3}).
 This completes
the proof of (\ref{g6.11}), and with it that of the lemma. \qed

Proposition \ref{p3.3} is a consequence of the next lemma.
\bl{l6.8} For all $n\in\Z$ we have
\begin{eqnarray*}
& &\Res_{z_{0}}z_{0}^{n}(z_{0}+z_{2})^{{\wt}a-1+{r\over T}+j}
\langle u',Y(a,z_{0}+z_{2})Y(b,z_{2})u\rangle\nonumber\\
&=&\Res_{z_{0}}z_{0}^{n}(z_{2}+z_{0})^{{\wt}a-1+{r\over T}+j}
\langle u',Y(Y(a,z_{0})b,z_{2})u\rangle.
\end{eqnarray*}
\el

\pf This is true for $n\geq -1$ by Lemma \ref{l6.7}. Let us write $n=-k+i$ with
$i\in \Z_+$ and proceed by induction $k.$ Induction yields
\begin{eqnarray*}
& &\Res_{z_{0}}z_{0}^{-k}(z_{0}+z_{2})^{{\wt}a-1+{r\over T}+j}
\langle u',Y(L(-1)a,z_{0}+z_{2})Y(b,z_{2})u\rangle\nonumber\\
&=&\Res_{z_{0}}z_{0}^{-k}(z_{2}+z_{0})^{{\wt}a-1+{r\over T}+j}
\langle u',Y(Y(L(-1)a,z_{0})b,z_{2})u\rangle.\label{g6.14}
\end{eqnarray*}

Using the residue property $\Res_{z}f'(z)g(z)\!+\!\Res_zf(z)g'(z)\!=\!0$
and the $L(-1)$-derivation property (\ref{2.5}) we have
\begin{eqnarray*}
& &\ \ \ \Res_{z_{0}}z_{0}^{-k}(z_{0}+z_{2})^{{\wt}a+{r\over T}+j}
\langle u',Y(L(-1)a,z_{0}+z_{2})Y(b,z_{2})u\rangle\nonumber\\
& &=-\Res_{z_{0}}\left({\partial\over\partial z_{0}}
z_{0}^{-k}(z_{0}+z_{2})^{{\wt}a+{r\over T}+j}\right)
\langle u',Y(a,z_{0}+z_{2})Y(b,z_{2})u\rangle\nonumber\\
& &=\Res_{z_{0}}k
z_{0}^{-k-1}(z_{0}+z_{2})^{{\wt}a+{r\over T}+j}
\langle u',Y(a,z_{0}+z_{2})Y(b,z_{2})u\rangle\nonumber\\
& &\ \ \ -\Res_{z_{0}}({\wt}a+{r\over T}+j)z_{0}^{-k}(z_{0}+z_{2})^{{\wt}a-1+
{r\over T}+j}
\langle u',Y(a,z_{0}+z_{2})Y(b,z_{2})u\rangle\nonumber\\
& &=\Res_{z_{0}}k
z_{0}^{-k-1}z_{2}(z_{0}+z_{2})^{{\wt}a-1+{r\over T}+j}
\langle u',Y(a,z_{0}+z_{2})Y(b,z_{2})u\rangle\nonumber\\
& &\ \ \ +\Res_{z_{0}}k
z_{0}^{-k}(z_{0}+z_{2})^{{\wt}a-1+{r\over T}+j}
\langle u',Y(a,z_{0}+z_{2})Y(b,z_{2})u\rangle\nonumber\\
& &\ \ \ -\Res_{z_{0}}({\wt}a+{r\over
T}+j)z_{0}^{-k}(z_{2}+z_{0})^{{\wt}a-1+{r\over T}+j}
\langle u',Y(Y(a,z_{0})b,z_{2})u\rangle\nonumber\\
& &=\Res_{z_{0}}k
z_{0}^{-k-1}z_{2}(z_{0}+z_{2})^{{\wt}a-1+{r\over T}+j}
\langle u',Y(a,z_{0}+z_{2})Y(b,z_{2})u\rangle\nonumber\\
& &\ \ \ +\Res_{z_{0}}k
z_{0}^{-k}(z_{2}+z_{0})^{{\wt}a-1+{r\over T}+j}
\langle u',Y(Y(a,z_{0})b,z_{2})u\rangle\nonumber\\
& &\ \ \ -\Res_{z_{0}}({\wt}a+{r\over
T}+j)z_{0}^{-k}(z_{2}+z_{0})^{{\wt}a-1+{r\over T}+j}
\langle u',Y(Y(a,z_{0})b,z_{2})u\rangle,\label{g6.15}
\end{eqnarray*}
and
\begin{eqnarray*}
& &\Res_{z_{0}}z_{0}^{-k}(z_{2}+z_{0})^{{\wt}a+{r\over T}+j}
\langle u',Y(Y(L(-1)a,z_{0})b,z_{2})u\rangle\nonumber\\
&=&-\Res_{z_{0}}\left({\partial\over\partial z_{0}}z_{0}^{-k}
(z_{2}+z_{0})^{{\wt}a+{r\over T}+j}\right)\langle
u',Y(Y(a,z_{0})b,z_{2})u\rangle
\nonumber\\
&=&\Res_{z_{0}}kz_{0}^{-k-1}(z_{2}+z_{0})^{{\wt}a+{r\over T}+j}
\langle u',Y(Y(a,z_{0})b,z_{2})u\rangle\nonumber\\
& &-\Res_{z_{0}}({\wt}a+{r\over T}+j)z_{0}^{-k}(z_{2}+z_{0})^{{\wt}a-1+{r\over
T}+j}
\langle u',Y(Y(a,z_{0})b,z_{2})u\rangle\nonumber\\
&=&\Res_{z_{0}}kz_{2}z_{0}^{-k-1}(z_{2}+z_{0})^{{\wt}a-1+{r\over T}+j}
\langle u',Y(Y(a,z_{0})b,z_{2})u\rangle\nonumber\\
& &+\Res_{z_{0}}kz_{0}^{-k}(z_{2}+z_{0})^{{\wt}a-1+{r\over T}+j}
\langle u',Y(Y(a,z_{0})b,z_{2})u\rangle\nonumber\\
& &-\Res_{z_{0}}({\wt}a+{r\over T}+j)z_{0}^{-k}(z_{2}+z_{0})^{{\wt}a-1+{r\over
T}+j}
\langle u',Y(Y(a,z_{0})b,z_{2})u\rangle.
\end{eqnarray*}
This yields the identity:
\begin{eqnarray*}
& &\Res_{z_{0}}z_{0}^{-k-1}(z_{0}+z_{2})^{{\wt}a-1+{r\over T}+j}
\langle u',Y(a,z_{0}+z_{2})Y(b,z_{2})u\rangle\nonumber\\
&=&\Res_{z_{0}}z_{0}^{-k-1}(z_{2}+z_{0})^{{\wt}a-1+{r\over T}+j}
\langle u',Y(Y(a,z_{0})b,z_{2})u\rangle,\label{3.15}
\end{eqnarray*}
and the lemma is proved.\qed

Let us now introduce an arbitrary ${1\over T}\Z_+$-graded $V[g]$-module
$M=\oplus_{n\in \1t\Z_+}M(n).$ As before we extend $M(0)^*$ to $M$ by letting
it annihilate $M(n)$ for $n>0.$ The proof of Proposition of \ref{p6.1}
with $\<u',\cdot\>$ suitably inserted gives:

\begin{prop}\label{p3.4}
Assume the following hold:

(i) $M=U(V[g])M(0).$

(ii) For $a\in V^r$ and  $u\in M(0)$
there is $k\in\Z$ such that
\begin{eqnarray}\label{ea4}
\langle u',(z_{0}+z_{2})^{k+{r\over T}}Y(a,z_{0}+z_{2})Y(b,z_{2})u\rangle
=\langle u',(z_{2}+z_{0})^{k+{r\over T}}Y(Y(a,z_{0})b,z_{2})u\rangle
\end{eqnarray}
for any $b\in V, u'\in M(0)^{*}$. Then in fact (\ref{ea4}) holds for any
$u\in M.$
\end{prop}

\begin{prop}\label{p3.5}
Let $M$ be as in Proposition \ref{p3.4}. Then for any $x\in
U(V[g]), a\in V^r, u\in M$, there is an integer $k$ such that
\begin{eqnarray}\label{ea5}
\langle u',(z_{0}\!+\!z_{2})^{k+{r\over T}}x\cdot
Y(a,z_{0}\!+\!z_{2})Y(b,z_{2})u\rangle\!=\!\langle
u',(z_{2}\!+\!z_{0})^{k+{r\over T}}x\cdot Y(Y(a,z_{0})b,z_{2})u\rangle
\end{eqnarray}
{\it for any $b\in V, u'\in M(0)^{*}$.}
\end{prop}

{\bf Proof.} Let $L$ be the subspace of $U(V[g])$ consisting of those
$x$ for which  (\ref{ea5}) holds. Let $x\in L$, let $c$ be any homogeneous
element of $V,$ and let $n\in {1\over T}{\Z}$. Then from (\ref{3.6'}) we have
\begin{eqnarray}
& &\ \ \ \langle u',xc_{n}Y(a,z_{0}+z_{2})Y(b,z_{2})u\rangle
(z_{0}+z_{2})^{k+{r\over T}}
\nonumber\\
& &=\sum_{i=0}^{\infty}\left(\begin{array}{c}n\\i\end{array}\right)
(z_{0}+z_{2})^{k+{r\over T}+n-i}\langle
u',xY(c_{i}a,z_{0}+z_{2})Y(b,z_{2})u\rangle
\nonumber\\
& &\ \ \ \ +\sum_{i=0}^{\infty}\left(\begin{array}{c}n\\i\end{array}\right)
z_{2}^{n-i}(z_{0}+z_{2})^{k+{r\over T}}\langle
u',xY(a,z_{0}+z_{2})Y(c_{i}b,z_{2})u\rangle\nonumber\\
& &\ \ \ \ +(z_{0}+z_{2})^{k+{r\over T}}\langle
u',xY(a,z_{0}+z_{2})Y(b,z_{2})c_{n}u\rangle.
\end{eqnarray}
The same method that was used in the proof of Proposition \ref{p3.4} shows
that
$xc_{n}\in L$. Since $U(V[g])$ is generated by all such
$c_{n}$'s, and since
(\ref{ea5}) holds for $x=1$ by Proposition \ref{p3.4},
we conclude that $L=U(V[g]),$ as desired. \qed

We can now finish the proof of  Theorem \ref{t6.1}.
We can take $M=M(U)$ in Proposition
\ref{p3.5}, as we may since $M(U)$ certainly satisfies the conditions placed
on $M$ prior to Proposition \ref{p3.4} and in Proposition \ref{p3.4}. Then
from the definition of $W$ (\ref{g6.3}), Lemma \ref{l6.4}
and Propositions \ref{p3.3}, \ref{p3.4} and \ref{p3.5} we conclude
that $U(V[g])W\subset J.$ So $\bar M(U)(0)\cong U.$

Turning to the proof of Theorem \ref{t6.3}, we have already seen that
$U(V[g])W\subset J.$ Then from Theorem \ref{t6.1} it is clear that
$L(U)=M(U)/J$ is a quotient of $\bar M(U)$ and hence an admissible
$g$-twisted $V$-module satisfying $L(U)(0)\cong U.$ Now clearly
$\O(L(U))\supset U,$ and if this is not an equality then there is
$n>0$ with $\O_n=\Omega(L(U))\cap L(U)(n)\ne 0.$ Then $0\ne U(V[g])\O_n$
intersects $U$ trivially, contradiction. So in fact $\O(L(U))\cong U.$
The remainder of Theorem \ref{t6.3} is straightforward to prove.
\qed

\section{Bijection between simple objects}

At this point we have a pair of functors $\O,L$ defined on appropriate
module categories:
\begin{equation}\label{7.1}
A_g(V)-\Mod
\begin{array}{l}
\stackrel{L}{\longrightarrow} \\
\stackrel{\Omega}{\longleftarrow}
\end{array}
\mbox{Adm}-g-V-\Mod
\end{equation}
Although $\O\circ L$ is equivalent to the identity, one cannot expect that
$L\circ\O$ is also equivalent to the identity in general. This is essentially
because there are examples of vertex operator algebras $V$ for which
the category of admissible $V$-modules contains objects which are not
completely reducible. For an example see [FZ]. We prove

\bl{l7.1} Suppose that $U$ is a simple $A_{g}(V)$-module.
Then $L(U)$ is a simple admissible $g$-twisted $V$-module.
\el

\pf If $0\ne W\subset L(U)$ is an admissible $g$-twisted
submodule then, by the definition of $L(U)$, we have
$W(0)=W\cap L(U)(0)\ne 0$. As $W(0)$ is an
$A_{g}(V)$-submodule of $U=L(U)(0)$ by Theorem \ref{t5.3} then
$U=W(0),$ whence $W\supset U(V[g])W(0)=U(V[g])U=L(U).$ \qed

\bt{t7.2} $L$ and $\O$ are equivalences when restricted to the full
subcategories of completely reducible $A_g(V)$-modules and
completely reducible admissible $g$-twisted $V$-modules respectively.
In particular, $L$ and $\Omega$ induces mutually
inverse bijections on the isomorphism classes of simple objects
in the category of $A_{g}(V)$-modules and admissible $g$-twisted $V$-modules
respectively.
\end{thm}

\pf We have $\Omega(L(U))\cong U$ for any $A_g(V)$-module by Theorem
\ref{t6.3}.

If $M$ is a completely reducible admissible
$g$-twisted $V$-module we must show $L(\O(M))\cong M.$ For this we may take
$M$ simple, whence $\O(M)$ is simple by Proposition
\ref{l2.9} (ii) and then $L(\O(M))$ is simple by Lemma \ref{l7.1}. Since both
$M$ and $L(\Omega(M))$ are simple quotients
of the universal object $\bar{M}(\Omega(M))$ then they are isomorphic
by Theorems \ref{t6.1} and \ref{t6.3}. \qed

\section{$g$-rational vertex operator algebras}

The definition of $g$-rational vertex operator algebra prior to Lemma
\ref{l3.4} says precisely that every object in the category of
admissible $g$-twisted $V$-modules is completely reducible. We have
the following amnibus result.
\bt{t8.1} Suppose that $V$ is a $g$-rational vertex operator algebra. Then the
following hold:

(a) $A_g(V)$ is a finite-dimensional, semi-simple associative algebra
(possibly 0).

(b) $V$ has only finitely many isomorphism classes of simple admissible
 $g$-twisted modules.

(c) Every simple admissible $g$-twisted $V$-module is an ordinary $g$-twisted
$V$-module.

(d) $V$ is $g^{-1}$-rational.

(e) The functors $L,\O$ are mutually inverse categorical equivalences
between the category of $A_g(V)$-modules and the category of admissible
$g$-twisted $V$-modules.

(f) The functors $L,\O$ induce mutually inverse categorical equivalences
between the category of finite-dimensional
$A_g(V)$-modules and the category of ordinary $g$-twisted $V$-modules.
\et

Suppose that (a) holds. Then all objects in the category of $A_g(V)$-modules
are completely reducible. Then (e) follows from Theorem \ref{t7.2}. Moreover
as $A_g(V)$ is of finite dimension it has only finitely many simple
modules, whence (b) follows from (e). Similarly (f) follows from (c).
So we must prove parts (a), (c), (d).

{\bf Proof of (a):} It suffices sto show that any
$A_g(V)$-module $U$ is completely reducible.

Now $L(U)$ is admissible and hence a direct sum of simple admissible
$g$-twisted $V$-modules. Application of the functor $\O$ shows that
$\O(L(U))$ is also completely reducible, so we are done since
$\O(L(U))\cong U.$

{\bf Proof of (c):} Let $M$ be a simple admissible $g$-twisted $V$-module. Then
$\O(M)$ is a simple $A_g(V)$-module, call it $U,$ and $L(U)\cong M.$

Now by Theorem \ref{t2.4} (iii), $\o+O_g(V)\in A_g(V)$ acts as
a scalar $h$, say, on $U.$ From the construction of $L(U)$ in Section
6 it follows that the graded subspaces $M(n)$ of $M$ are precisely
the distinct eigenspaces of $L(0)$ on $M.$ That is, $L(0)$ is semi-simple
as an operator on $M$ and for $n\in\1t\Z_+$ we have
\be{8.1}
M(n)=\{m\in M|L(0)=(n+h)m\}.
\end{equation}

Next, if we think of $U$ as a left $A_g(V)$-module then $U^*$ is naturally
a left $A_g(V)^{opp}$-module (opposite algebra). Theorem \ref{t2.4}
(ii) tells us that there is a canonical algebra isomorphism $A_{g^{-1}}(V)
\cong  A_g(V)^{opp},$
so $U^*$ is a left $A_{g^{-1}}(V)$-module. It is simple by part (a). Now
apply Theorem \ref{t6.3} (with $g^{-1}$ in place of $g$): we conclude
that $L(U^*)$ is an admissible $g^{-1}$-twisted $V$-module. Moreover
$L(U^*)$ is simple by Lemma \ref{l7.1}.

Let $N=L(U^*)'.$ So $N$ is an admissible $g$-twisted $V$-module (see
Lemma \ref{l3.8}). We will show that $N$ is also simple.

Let $W$ be the admissible $g$-twisted submodule of $N$ generated by $N(0)
=(U^*)^*=U.$ As $V$ is $g$-rational then $N=W\oplus W_0$ for some
admissible $g$-twisted submodule $W_0$ of $N.$ Obviously $W_0\cap U=0,$
so $\<U^*,W_0\>=0.$ As $U^*$ generates $L(U^*)$ we get $W_0=0.$
So $N=W$ is generated by $U.$

Now if $X$ is any non-zero submodule of $N$ then
we have $X\cap U\ne 0$ as $N$ is completely reducible. Since $U$ is a
simple
$A_g(V)$-module then $U\subset X,$ whence $U=X.$ So indeed $N$ is simple.
Then both $N$ and $L(U)$ are simple admissible $g$-twisted $V$-modules
generated by $U,$ so we must have $N\cong L(U)$ (cf. Theorems \ref{t6.1}
and \ref{t6.3}).

Applying Lemma \ref{l3.7} shows that $N$ has countable dimension, so
the same is true of each graded subspace. Thus $(L(U^*)(n))^*$ has
countable dimension.  This can only happen if $L(U^*)(n)$ is of {\em
finite} dimension.

We now deduce that in fact $L(U^*)$ is an ordinary $g^{-1}$-twisted $V$-module.
Indeed (\ref{8.1}) applies to $L(U^*)$ also, so that axiom (\ref{g3.11})
is fulfilled. Then as the contragredient module of $L(U^*),$ $L(U)$ is also
an ordinary $g$-twisted $V$-module. This completes the proof of part
(c).

{\bf Proof of (d):} Let $\{W^{1}, \cdots, W^{k}\}$ be representatives
for the equivalence
classes of simple $g$-twisted $V$-modules

Consider any ${1\over T}{\Z}_{+}$-graded weak $g^{-1}$-twisted $V$-module
$M=\oplus_{n\in {1\over T}{\Z}_{+}}M(n)$.
Then $M'=\oplus_{n\in {1\over T}{\Z}_{+}}M(n)^{*}$ is
an admissible $g$-twisted $V$-module which
is completely reducible. We can write
$$M'\simeq U^{1}\otimes W^{1}\oplus U^{2}\otimes W^{2}\oplus \cdots \oplus
U^{k}\otimes W^{k}$$
for certain vector spaces $U^i.$ Clearly $U^i$ are ${1\over T}\Z_+$-graded:
$$U^i=\displaystyle{\oplus_{m\in{1\over T}\Z_+}U^i(m)}$$
such that
$$M(n)^*=\displaystyle{\oplus_{i=1}^k\oplus_{s,t\in {1\over T}\Z_+,
s+t=n}U^i(s)\otimes W^i(t)}.$$
As $W^i$ is an ordinary $g$-twisted
module, each $W^i(t)$ is
finite-dimensional. It follows that
$$(M(n)^*)^*=\displaystyle{\oplus_{i=1}^k\oplus_{s,t\in {1\over T}\Z_+,
s+t=n}U^i(s)^*\otimes W^i(t)^*}$$
and thus
$$(M')'=\displaystyle{\oplus_{i=1}^k\oplus_{n\in {1\over T}\Z_+}U^i(n)^*\otimes
(W^i)'.}$$
So $(M')'$ is completely reducible, whence so to is $M\subset (M')'.$
This completes the proof of Theorem \ref{t8.1}. \qed

\section{Further applications}

It is a well-known conjecture that $V$ always possess at least one
ordinary $g$-twisted module. (Here $V$ is any vertex operator algebra
and $g$ an automorphism of order $T.$) Somewhat weaker is the
conjecture that $A_g(V)$ is  non-zero; this is equivalent to the
existence of a simple admissible $g$-twisted $V$-module by Theorem
\ref{t7.2}. We have the following contribution to this problem:

\bt{t9.1} Suppose that $A_{g}(V)$ is of finite dimension. Then there
is at least one simple $g$-twisted $V$-module.
\et

We begin the proof with a variation on the theme of Section 8.
\bl{l9.2}
Let $M=\oplus_{n\in\1t\Z_+}M(n)$ be an admissible $g$-twisted $V$-module
such that (\ref{8.1}) holds, i.e., $M(n)$ is the $n+h$-eigenspace
for the action of $L(0).$ Assume that the contragredient
module $M'$ is simple. Then $M$ is a simple (ordinary) $g$-twisted module.
\el

\pf  If $W$ is  a non-zero submodule of $M$ then
$W^{\perp}=\{f\in M'| \<f,W\>=0\}$
is a weak $g^{-1}$-twisted submodule of $M'$.
As $M'$ is simple then $W^{\perp}=0,$ so $W=M.$ This shows that $M$ is
simple.

As $M'$ is simple then it has countable dimension by Lemma \ref{l3.7}, and the
remainder of the proof now follows
as in the proof of Theorem \ref{t8.1} (c). \qed

Denote by $\M_g(V)$ the (equivalence classes of) simple admissible
$g$-twisted $V$-modules. It is a finite set since we are now assuming
that $A_g(V)$ is of finite dimension. As before, $L(0)$ is semi-simple
as an operator on these modules, and we denote by $S_{g}(V)\subset\C$
the set of {\em lowest weights}, i.e., the set of eigenvalues $h$ in the
notation of (\ref{8.1}). Define a partial order $\geq $ on $S_g(V)$
as follows: $h_{1}\ge h_{2}$ if, and only if, $h_{1}-h_{2}\in {1\over
T}{\Z}_{+}$. Let $S_g^*(V)$ be the maximal elements in the partial order, and
$\M^*_g(V)$ the modules of $\M_g(V)$ whose lowest weights lie in $S_g^*(V).$
\bl{l9.3} Suppose that $M\in \M_g^*(V).$ Then $M$ is a $g$-twisted $V$-module.
\el

\pf Since $M$ is simple then $\O(M)=M(0)$ is a simple $A_g(V)$-module, hence of
finite dimension. Then $M(0)^*$ is a simple $A_{g^{-1}}(V)$-module
(cf. the proof of part (c) of Theorem \ref{t8.1}) and $M(0)^*=M'(0).$

Now as $A_{g^{-1}}\cong A_g(V)^{opp}$ (Theorem \ref{t2.4} (ii)) then
clearly $S_{g}(V)$ and $S_{g^{-1}}(V)$ are equal as sets. So we must be in the
situation that $M'(0)$ generates $M'.$ As usual this leads to the
conclusion that $M'$ is simple. Now the present lemma follows from Lemma
\ref{l9.2}. \qed

Notice that Theorem \ref{t9.1} is a special case of Lemma \ref{l9.3}.

In the following, for any type of $V$-module for which $L(0)$ is semi-simple,
we let $M_h$ denote the $h$-eigenspace of $L(0).$

Next  we give a sufficient condition for the existence of a composition series
of finite length for any $g$-twisted $V$-module.

\begin{prop}\label{p3.10}
Suppose that there are only finitely many simple $g$-twisted
$V$-modules (up to equivalence). Then any $g$-twisted $V$-module $M$ has
a finite composition series such that each factor is simple.
\end{prop}

{\bf Proof.} Let $\{W^{1}, \cdots, W^{k}\}$ be the set of equivalence
classes of simple $g$-twisted $V$-modules and let $h_{i}$ be the
lowest weight of $W^{i}.$ Let $M$ be a $g$-twisted $V$-module.

{\bf Claim 1:} {\em If $h$ is a minimal weight of $M$ in the sense
that $M_h\ne 0$ and $M_{h-n}=0$ for all
positive $n\in \frac{1}{T},$ then $h\in \{h_{1},\cdots,
h_{k}\}$.} Suppose $h$ is a minimal weight of $M$. Then $M_{h}\subseteq
\Omega(M)$ and $M_{h}$ is an $A_{g}(V)$-submodule of $\Omega(M)$. Let
$U$ be an irreducible $A_{g}(V)$-submodule of $M_{h}$ and let $W$ be
the $g$-twisted $V$-submodule generated by $U$. Then $W$ has a
unique simple quotient module, and its  lowest weight is $h.$
 Thus $h=h_{i}$ for some $i$.

{\bf Claim 2:} {\em $M$ is generated by
$E(M)=\oplus_{i=1}^{k} M_{h_{i}}$ (note that $E(M)$ is
finite-dimensional).}  Let $W$ be the $g$-twisted $V$-submodule
generated by $E(M)$. If $M\ne W$, $M/W$ is a (nonzero) $g$-twisted
$V$-module such that no homogeneous subspace has weight $h_{i}$ for any $i.$
So there must be a lowest weight $h$ of $M/W$ such that $h\ne
h_{i}$ for any $1\le i\le k$. Now apply Claim 1 to get a contradiction.

{\bf Claim 3:} {\em If $W^{1}$ and $W^{2}$ are two submodules of
$M$ such that $W^{1}\cap E(M)=W^{2}\cap E(M)$,
then $W^{1}=W^{2}$.} Since $E(W^{1})=W^{1}\cap E(M)=W^{2}\cap
E(M)=E(W^{2})$, this follows from Claim 2 immediately.

Let $\cal{S}$ be the set of all submodules of $M$ partially ordered
by inclusion.

{\bf Claim 4:} {\em There exists a finite maximal chain in $\cal{S}$.}
It follows from Zorn's Lemma that there exist maximal
chains. Let $\cdots \subseteq M^{-1}\subseteq M^{0}\subseteq
M^{1}\subseteq \cdots$ be an ascending chain in
$\cal{S}$. Then we have:
$$\cdots\subseteq E(M)\cap M^{-1}\subseteq
E(M)\cap M^{0}\subseteq E(M)\cap M^{1}\subseteq \cdots$$ Since
$E(M)$ is finite-dimensional, there are nonnegative integers $m$ and
$n$ such that
$$E(M)\cap M^{s}=E(M)\cap M^{-m},\;E(M)\cap
M^{n}=E(M)\cap M^{t}$$ for any $s\le-m$ and $t\ge n.$ Thus by Claim
3, $M^{s}=M^{-m}$ and $M^{n}= M^{t}.$

It is clear that for any maximal chain, all the factors are
simple $g$-twisted $V$-modules. The proof is
complete. \qed

A vertex operator algebra $V$ is said to be {\em holomorphic} if $V$
is the only simple $V$-module up to equivalence. The famous
moonshine module vertex operator algebra $V^{\natural}$ [FLM] is
an example [D2]. The following proposition is an application of Proposition
\ref{p3.10}.

\begin{prop}\label{p3.11}
Suppose that $V$ is a holomorphic VOA and that $V$ contains a
rational vertex operator subalgebra $U$ (with the same Virasoro
element). Then any $V$-module is completely reducible.
\end{prop}

We need the following Lemma from [L1].
\bl{lv} Let $V$ be a vertex operator algebra and $M$ a weak
$V$-module.

(i) Let $m\in M$ be a vacuum-like vector, that is, $L(-1)m=0.$
Then the weak submodule of $M$ generated by $m$ is isomorphic
to $V,$ the isomorphism arising via the map $m\mapsto {\bold 1}.$

(ii) Conversely, if $V$ is isomorphic to $M$ as weak $V$-modules, then
the image of ${\bold 1}$ in $M$ is a vacuum-like vector.
\el

{\bf Proof of Proposition \ref{p3.11}:} Let $M$ be a $V$-module. By
 Proposition \ref{p3.10}, there is a composition series $0\subseteq
 M^{1}\subseteq M^{2}\subseteq \cdots \subseteq M^{k}=M$.  Assume that
 $M^{i}$ is completely reducible for some $i.$ {}From the assumption,
 there is a $U$-submodule $B^{i}$ of $M^{i+1}$ such that
 $M^{i+1}=B^{i}\oplus M^{i}$.  Since $M^{i+1}/M^{i}\simeq V$, there is
 $0\ne u\in M^{i+1}$ such that $u+M^i$ is a vacuum-like vector in
 $M^{i+1}/M^i$ by Lemma \ref{lv}.  Let $u=a+b$ where $a\in B^{i}, b\in
 M^{i}$. Then $L(-1)u=L(-1)a+L(-1)b\in M^i.$ Thus $a\ne 0$ and
 $L(-1)a=0$ as $L(-1)$ preserves both $M^i$ and $B^i.$ This show that
 $a$ is a vacuum-like vector in $M^{i+1}.$ Again by Lemma \ref{lv},
 $a$ generates a $V$-submodule of $M^{i+1}$ isomorphic to $V$. Then
 $M^{i+1}$ is the sum of $M^{i}$ and $U(V[1])a$.  So $M^{i+1}$ is
 completely reducible. It follows by induction that each $M^{i}$ is
 completely reducible for $i\ge 1$. In particular, $M$ is completely
 reducible.
\qed

\br{r4.1} The moonshine module vertex operator algebra $V^{\natural}$ [FLM]
has a rational vertex operator subalgebra
$U$ which satisfies the conditions of Proposition \ref{p3.11}.
It is isomorphic to the tensor product of 48 irreducible highest weight
modules
for the Virasoro algebra with central charge $1\over 2$ [DMZ]. Proposition
\ref{p3.11} thus gives a proof that any $V^{\natural}$-module is
completely reducible which is shorter than the original proof [D2].
\er

It is conjectured in [FLM] that any holomorphic vertex operator algebra $V$
of rank 24 with $V_1=0$ is isomorphic to the moonshine module vertex operator
algebra $V^{\natural}$ in [FLM]. The following proposition asserts that
any ordinary module for such a vertex operator
algebra is completely reducible. When applied to $V^{\natural}$ itself,
this
gives another proof of complete reducibility of any $V^{\natural}$-module.

\begin{prop}\label{p3.11'}
Suppose that $V$ is a holomorphic VOA such that $V_1=0.$ Then any $V$-module
is completely reducible.
\end{prop}

{\bf Proof.} Let $M$ be any $V$-module and $W\subset M_0$ be the subspace
of vacuum-like vectors (cf. Lemma \ref{lv}). Let $M'$ be the submodule of
$M$ generated
by $W.$ Then $M'$ is completely reducible. If $M'\ne M$ consider $M/M'.$
Let $u\in M\setminus M'$ such that $u+M'$ is a vacuum-like vector, that
is $L(-1)u\in M'.$ Note that $u\in M_0$ and $L(-1)u\in M_1.$ Since
$M$ has a finite composition series and each factor is isomorphic to $V$
we see that $M_1=0.$ Thus $L(-1)u=0$ and
$u$ is a vacuum-like vector. This is a contradiction because $u$ is not
in $W.$ \qed

\begin{prop}\label{p3.23}
Suppose that $V^{0}$ contains a rational vertex operator
subalgebra $U$ (with the same Virasoro element) such that
the fusion rules among any three irreducible
$U$-modules is finite. Then any simple admissible $g$-twisted $V$-module
is an ordinary
$g$-twisted $V$-module.
\end{prop}

\pf Let $M$ be a simple admissible $g$-twisted $V$-module with lowest weight
$h$.
Since $M$ is a
completely reducible $U$-module we can take a simple admissible
$U$-submodule $W$ of $M$.  Then by Proposition 2.4 of [DM2]
or Lemma 6.1.1 of [L2] $M$ is linearly
spanned by the coefficients of $Y_{M}(a,z)u$ for $a\in V$ and fixed $u\in W.$
Regarding $V$, $W$ and $M$ as $U$-modules, we have an intertwining operator
$Y_M$ of type ${M\choose V\,W}$ (see [FHL] for the definition of intertwining
operator). It follows from the universal property of the  tensor product
that there is a $U$-homomorphism $\psi$ from $V\boxtimes W$ onto $M$ (cf.
[HL0]-[HL1] and [L2]). From our assumption,
$V\boxtimes W$ is
a sum of finitely many irreducible $U$-modules, so that any homogeneous
subspace is finite-dimensional.
Then any homogeneous subspace of $M$ is finite-dimensional. That is, $M$ is an
ordinary $g$-twisted
$V$-module. \qed

A similar result has been obtained in [H] in the special case when $g=1$ and
$V$ contains a rational vertex operator subalgebra which
is a tensor product  of vertex operator algebras associated with
the highest weight irreducible representations for the discrete series
of the Virasoro algebra.

\begin{prop}\label{p3.24}
Suppose that $V$ is a holomorphic VOA and that $V$ contains a rational
vertex operator
subalgebra $U$ (with the same Virasoro element)
such that the fusion rules among any three irreducible
$U$-modules are finite. Then $V$ is rational.
\end{prop}

{\bf Proof.} We need to prove that any admissible
$V$-module $M$ is completely reducible. Now both $V$ and $M$ are
direct sums of simple $U$-modules as $U$ is rational.
Let $W$ be a
simple $U$-submodule of $M$ and let $\bar{W}$ be the weak
$V$-submodule generated by $W$.  As in the proof of
Proposition \ref{p3.23}, we easily
show that any homogeneous subspace of the
$U$-module $V\boxtimes W$ is finite-dimensional. Thus being a
$U$-homomorphic image of a $U$-module $V\boxtimes W$, $\bar{W}$ is a
$U$-module.  By Proposition \ref{p3.11}, $\bar{W}$ is a completely
reducible $V$-module. Thus $M$ is a completely reducible
$V$-module. \qed

\begin{rem}\label{r3.25} (i) Recall Remark \ref{r4.1}. The moonshine module
vertex operator algebra $V^{\natural}$ satisfies the condition of
Proposition \ref{p3.24}. So $V^{\natural}$ is a rational holomorphic
vertex operator algebra.

(ii) A situation in which Proposition \ref{p3.23} applies is where
$V=V^{\natural}$ and $g$ any involution in the Monster. This is studied
in detail in [DLM2].
\end{rem}


\begin{thebibliography}{FLM1}


\bibitem[B]{2}
R. E. Borcherds, Vertex algebras, Kac-Moody algebras, and the Monster,
{\it Proc. Natl. Acad. Sci. USA} {\bf 83} (1986), 3068-3071.

\bibitem[DHVW]{1} L. Dixon, J. Harvey, C. Vafa and E. Witten,
Strings on orbifolds, {\em Nucl. Phys.} {\bf B261} (1985),651; II,
{\em Nucl. Phys.} {\bf B274} (1986),285.

\bibitem[DGM]{1} L. Dolan, P. Goddard and P. Montague, Conformal field
theory of twisted vertex operators, {\em Nucl. Phys.} {\bf B338} (1990),
529.

\bibitem[D1]{}
C. Dong, Twisted modules for vertex algebras associated with even lattices,
{\it J. of
Algebra} {\bf 165} (1994), 91-112.

\bibitem[D2]{1}
C. Dong, Representations of the moonshine module vertex operator
algebra, {\it Contemporary Math.},  {\bf 175} (1994),
27-36.

\bibitem[DL]{1}
C. Dong and J. Lepowsky, Generalized Vertex Algebras and Relative
Vertex Operators, Progress in Math., {\bf Vol. 112}, Birkhauser,
Boston, 1993.

\bibitem[DLM1]{1} C. Dong, H. Li and G. Mason, Regularity of rational
vertex operator algebras, preprint, q-alg/9508018.

\bibitem[DLM2]{2} C. Dong, H. Li and G. Mason, Some twisted sectors for
the Moonshine Module, {\it Contemporary Math.}, to appear, q-alg/9504014.

\bibitem[DM1]{1}
C. Dong and G. Mason, On quantum Galois theory,
preprint, 1994; hep-th/9412037.

\bibitem[DM2]{2}
C. Dong and G. Mason, On the operator content of nilpotent
orbifolds, preprint, 1994; hep-th/9412109.

\bibitem[DMZ]{1}
C. Dong, G. Mason and Y. Zhu, Discrete series of the Virasoro
algebra and the moonshine module,  {\em Proc. Symp. Pure. Math., American Math.
Soc.} {\bf 56} II (1994), 295-316.

\bibitem[FFR]{1}
Alex J. Feingold, Igor B. Frenkel and John F. X. Ries, {\it Spinor
Construction of Vertex Operator Algebras, Triality, and
$E_{8}^{(1)}$}, {\em Contemporary Math.} {\bf 121} (1991).

\bibitem[FHL]{1}
I. Frenkel, Y.-Z. Huang and J. Lepowsky, On axiomatic approaches to
vertex operator algebras and modules, Memoirs Amer. Math.
Soc. {\bf 104}, 1993.

\bibitem[FLM]{1}
I. Frenkel, J. Lepowsky and A. Meurman, {\it Vertex Operator Algebras
and the Monster}, Pure and Appl. Math., {\bf Vol. 134}, Academic Press,
Boston, 1988.

\bibitem[FZ]{1}
I. Frenkel and Y. Zhu, Vertex operator algebras associated to
representations of affine and
Virasoro algebras, {\it Duke Math. J.} {\bf 66} (1992), 123-168.


\bibitem[H]{1}
Y. Huang, Virasoro vertex operator algebras, the (nonmeromorphic)
operator product expansion and the tensor product theory, preprint.

\bibitem[HL0]{1}
Y. Huang and J. Lepowsky, Toward a theory of tensor product for
representations
for a vertex operator algebra, in {\it Proc. 20th
International Conference on Differential Geometric Methods in
Theoretical Physics, New York, 1991,} ed. S. Catto and A. Rocha, World
Scientific, Singapore, 1992, {\bf Vol. 1}, 344-354.

\bibitem[HL1]{1}
Y. Huang and J. Lepowsky, A theory of tensor product for
module category of a vertex operator algebra, I,II, preprint (1993).

\bibitem[L1]{1}
H. Li, Smmetric invariant bilinear forms on vertex
operator algebras, {\em Journal of Pure and Applied Algebra}, {\bf 96} (1994),
279-297.

\bibitem[L2]{1}
H. Li, Representation theory and tensor product theory for vertex  operator
algebras, Ph.D. thesis, Rutgers University, 1994.

\bibitem[MS]{MS} G. Moore and N. Seiberg, Classical and quantum
conformal field theory, {\em Comm. Math. Phys.} {\bf 123} (1989), 177-254.

\bibitem[T]{1} H. Tamanoi, Elliptic genera and vertex operator superalgebras,
preprint.

\bibitem[X]{1}
X. Xu, Twisted modules for colored vertex operator superalgebras,
preprint, 1992.

\bibitem[Z]{1}
Y. Zhu, Modular invariance of characters of vertex operator algebras,
{\em J. Amer, Math. Soc.,} to appear.
\end{thebibliography}
\end{document}